\documentclass[twoside,12pt]{article}
\usepackage{amssymb}
\usepackage{amsthm}
\usepackage{amsmath}
\usepackage{a4wide}
\usepackage[all]{xy}
\usepackage{fancyhdr}

\newtheorem{theorem}{Theorem}[section]
\newtheorem{proposition}[theorem]{Proposition}
\newtheorem{corollary}[theorem]{Corollary}
\newtheorem{lemma}[theorem]{Lemma}
\theoremstyle{definition}
\newtheorem*{notation}{Notation}
\newtheorem*{Beweis}{Proof}
\newtheorem{definition}[theorem]{Definition}
\newtheorem{punto}[theorem]{}
\theoremstyle{remark}
\newtheorem{remark}[theorem]{Remark}
\newtheorem{ex}[theorem]{Example}
\newtheorem{exs}[theorem]{Examples}
\newtheorem{c-ex}[theorem]{Counterexample}
\newtheorem{c-exs}[theorem]{Counterexamples}
\newtheorem{remarks}[theorem]{Remarks}
\CompileMatrices  \swapnumbers
\CompileMatrices
\input{tcilatex}
\begin{document}

\title{Zariski Topologies for Coprime and Second Submodules\thanks{\textbf{%
MSC2000:}\ 16D10 (13C05, 13C13, 54B99) \newline
\textbf{Key Words:\ }Zariski Topology; Dual Zariski Topology; Multiplication
Module; Comultiplication Module; Coprime Module; Second Module; Dual Ring;
Strongly Hollow Submodule; Strongly Irreducible Submodule}}
\author{\textbf{Jawad Abuhlail}\thanks{%
The author would like to acknowledge the support provided by the Deanship of
Scientific Research (DSR) at King Fahd University of Petroleum $\&$ Minerals
(KFUPM) for funding this work through project No. FT08009.} \\
Department of Mathematics and Statistics\\
King Fahd University of Petroleum $\&$ Minerals \\
abuhlail@kfupm.edu.sa}
\date{\today }
\maketitle

\begin{abstract}
Let $M$ be a non-zero module over an associative (not necessarily
commutative) ring. In this paper, we investigate the so-called \emph{second}
and \emph{coprime} submodules of $M.$ Moreover, we topologize the spectrum $%
\mathrm{Spec}^{\mathrm{s}}(M)$ of second submodules of $M$ and the spectrum $%
\mathrm{Spec}^{\mathrm{c}}(M)$ of coprime submodules of $M,$ study several
properties of these spaces and investigate their interplay with the
algebraic properties of $M.$
\end{abstract}

\section{Introduction}

\qquad Several papers considered the so called \emph{top modules,} i.e.
modules over commutative rings whose spectrum of \emph{prime submodules}
attains a Zariski-like topology, e.g. \cite{Lu1995, Lu1997, Lu1999, MMS1997,
MMS1998}. In \cite{Abu2006} and \cite{Abu2008}, the author investigated and
topologized the spectrum of \textit{fully coprime subbicomodules} of a given
non-zero duo bicomodule over a coring. Recently, he introduced module
theoretic versions of these results in \cite{Abu2011}, where a \textit{dual
Zariski topology }was introduced on the spectrum of \textit{fully coprime
submodules }of a given non-zero duo module over an associative ring.
Moreover, he introduced and studied a Zariski topology on the spectrum of
\textit{fully prime submodules} of a given non-zero duo module in \cite{Abu}.

As a dual notion of \textit{prime submodules,} Yassemi \cite{Yas2001}
introduced the notion of \textit{second submodules }of a given non-zero
module over a commutative ring. This notion was generalized to modules over
arbitrary associative rings by Annin in \cite{Ann2002}, where a \textit{%
second module} was called a \textit{coprime module}. Moreover, the notion of
\textit{coprime submodules} was introduced by Kazemifard et al. \cite{KNR}.
In this paper, we investigate conditions under which the spectrum $\mathrm{%
Spec}^{\mathrm{s}}(M)$ of second submodules ($\mathrm{Spec}^{\mathrm{c}}(M)$
of coprime submodules) of a given non-zero \emph{module} $M$ over an
arbitrary associative -- not necessarily commutative -- ring $R$ attains a
(dual) Zariski topology. We study these spaces and investigate the interplay
between the properties of these topologies and the algebraic properties of $%
_{R}M.$

After this introductory section, we introduce in Section 2 some
preliminaries. In particular, we recall some properties and notions of
modules that will be needed in the sequel. Section 3 is devoted to a study
of the second and coprime submodules of $M.$ In Section 4, we introduce a
\textit{dual Zariski topology} on $\mathrm{Spec}^{\mathrm{s}}(M)$ and study
its properties. The results obtained are similar to results on the spectrum $%
\mathrm{Spec}^{\mathrm{fc}}(M)$ of \textit{fully coprime submodules} of $%
_{R}M$ \cite{Abu2011}. In Section 5, we investigate a \textit{Zariski
topology} on $\mathrm{Spec}^{\mathrm{c}}(M).$ The results obtained are
similar to results on the spectrum $\mathrm{Spec}^{\mathrm{fp}}(M)$ of
\textit{fully prime submodules} of $_{R}M$ \cite{Abu2011}.

Throughout, $R$ is an associative (not necessarily commutative) ring with $%
1_{R}\neq 0_{R}.$ With $\mathrm{Max}(R)$ (resp. $\mathrm{Max}(_{R}R),$ $%
\mathrm{Max}(R_{R})$) we denote the spectrum of maximal ideals (resp.
maximal left ideals, maximal right ideals) of $R.$ On the other hand, we
denote by $\mathrm{Min}(R)$ (resp. $\mathrm{Min}(_{R}R),$ $\mathrm{Min}%
(R_{R})$) the set of minimal ideals (resp. minimal left ideals, minimal
right ideals) of $R.$ Recall that an ideal $\mathfrak{p}$ of $R$ is said to
be (\textit{completely}) \textit{prime} iff for any ideal $I,J$ of $R$ (any $%
a,b\in R$) with $IJ\subseteq \mathfrak{p}$ ($ab\in \mathfrak{p}$), either $%
I\subseteq \mathfrak{p}$ or $J\subseteq \mathfrak{p}$ ($a\in \mathfrak{p}$
or $b\in \mathfrak{p}$). With $\mathrm{Spec}(R)$ ($\mathrm{CSpec}(R)$) we
denote the spectrum of (completely) prime ideals of $R.$ With $\mathrm{Rad}%
(R):=\dbigcap\limits_{\mathfrak{m}\in \mathrm{Max}(R)}\mathfrak{m}$ ($%
\mathrm{Prad}(R):=\dbigcap\limits_{\mathfrak{p}\in \mathrm{Spec}(R)}%
\mathfrak{p}$) we denote the (\textit{prime})\textit{\ radical} of $R.$ The
set of zero-divisors in $R$ is denoted by $Z(R)$ while the group of
invertible elements of $R$ is denoted by $U(R).$ Unless otherwise explicitly
mentioned, a module will mean a \textit{left} $R$-module and an ideal is a
\textit{two-sided} ideal.

Moreover, we fix an arbitrary left non-zero $R$-module $_{R}M$ with ring of
endomorphisms $S:=\mathrm{End}(_{R}M)^{op}$ and consider $M$ as an $(R,S)$%
-bimodule in the canonical way. We write $K\leq _{R}M$ ($K\lvertneqq _{R}M$)
to indicate that $L$ is a (proper) submodule of $M$ and denote with $\pi
_{K}:M\rightarrow M/K$ the canonical surjection. With $\mathbb{P}\subset
\mathbb{Z}$ we denote the set of prime positive integers.

\section{Preliminaries}

\qquad For the convenience of the reader, we recall in this section some
definitions and properties of modules that will be used in the sequel.
Moreover, we illustrate these notions by introducing several examples. For
more information, the interested reader may refer to any book in Module
Theory (e.g. \cite{Wis1991}).

\begin{punto}
We call $L\leq _{R}M$ \textit{fully invariant} iff $L$ is also an $S$%
-submodule. We call $_{R}M$ \textit{duo} iff every $R$-submodule of $M$ is
fully invariant. The ring $R$ is said to be \textit{left duo} (\textit{right
duo}) iff every left (right) ideal is two-sided and to be \textit{left
quasi-duo} (\textit{right quasi-duo}) iff every maximal left (right) ideal
of $R$ is two-sided. Moreover, $R$ is said to be (\textit{quasi-}) \textit{%
duo} iff $R$ is left and right (quasi-) duo.
\end{punto}

\begin{notation}
With $\mathcal{L}(M)$ ($\mathcal{L}^{\mathrm{f.i.}}(M)$) we denote the
lattice of (fully invariant) $R$-submodules of $M$ and with $\mathcal{I}%
_{r}(R)$ (resp. $\mathcal{I}_{l}(R),$ $\mathcal{I}(R)$) the lattice of right
(resp. left, two-sided) ideals. For subsets $X,Y\subseteq M$ and $Z\subseteq
R$ we set%
\begin{equation*}
\begin{tabular}{lllll}
$(X$ & $:_{Z}$ & $Y)$ & $:=$ & $\{r\in Z\mid ry\in X\text{ for every }y\in
Y\};$ \\
$(X$ & $:_{Y}$ & $Z)$ & $:=$ & $\{m\in Y\mid rm\in X\text{ for every }r\in
Z\}.$%
\end{tabular}%
\end{equation*}%
In particular, $\mathrm{ann}_{R}(Y):=(0:_{R}Y)$ and $\mathrm{ann}%
_{M}(Z):=(0:_{M}Z).$ On the other hand, for any non-empty subsets $%
K\subseteq M$ and $I\subseteq $ $S$ we set $\mathrm{An}(K):=(0:_{S}K)$ and $%
\mathrm{Ke}(I):=(0:_{M}I).$ Moreover, we set%
\begin{equation*}
\mathcal{L}_{m}(M):=\{IM\mid I\in \mathcal{I}(R)\}\text{ and }\mathcal{L}%
_{c}(M):=\{L\leq _{R}M\mid L=(0:_{M}(0:_{R}L))\}.
\end{equation*}
\end{notation}

\begin{definition}
\label{injectivity}We say $_{R}M$ is

\textit{self-injective} iff for every $K\leq _{R}M,$ every $f\in \mathrm{Hom}%
_{R}(K,M)$ extends to some $\widetilde{f}\in S;$

\textit{intrinsically injective} iff $\mathrm{AnKe}(I)=I$ for every finitely
generated right ideal $I$ of $S;$

\textit{self-cogenerator} iff $M$ cogenerates all its factor $R$-modules.
\end{definition}

\begin{punto}
By $\mathrm{Max}(M)$ ($\mathrm{Max}^{\mathrm{f.i.}}(M)$), we denote the
\textit{possibly empty} class of maximal $R$-submodules of $M$ (the class of
maximal $(R,S)$-subbimodules of $_{R}M_{S}$). For every $L\leq _{R}M,$ we set%
\begin{equation*}
\mathcal{M}(L):=\{K\in \mathrm{Max}(M)\mid K\supseteq L\}\text{ and }%
\mathcal{M}^{_{\mathrm{f.i.}}}(L):=\{K\in \mathrm{Max}^{\mathrm{f.i.}%
}(M)\mid K\supseteq L\}.
\end{equation*}
\end{punto}

\begin{punto}
Let $L\leq _{R}M.$ We say that $L$ is \emph{superfluous} or \emph{small} in $%
M,$ and write $L\ll M,$ iff $L+\widetilde{L}\neq M$ for every $\widetilde{L}%
\lvertneqq _{R}M.$ The \emph{radical} of $M$ is defined as
\begin{equation*}
\mathrm{Rad}(M):=\dbigcap\limits_{L\in \mathrm{Max}(M)}L=\dsum\limits_{L\ll
M}L\text{ \ (}:=M\text{ iff }\mathrm{Max}(M)=\varnothing \text{).}
\end{equation*}
\end{punto}

\begin{definition}
We say $_{R}M$ is

\emph{local} iff $M$ contains a proper $R$-submodule that contains every
proper $R$-submodule of $M,$ equivalently iff $\dsum\limits_{L\lvertneqq
_{R}M}L\neq M$ (this is also equivalent to $_{R}M$ being cyclic, or finitely
generated, and having a unique maximal submodule);

\emph{hollow} (or \emph{couniform}) iff for any $L_{1},L_{2}\lneqq _{R}M$ we
have $L_{1}+L_{2}\lvertneqq _{R}M$, equivalently iff every proper $R$%
-submodule of $M$ is superfluous;

\emph{coatomic} (or $B$\textit{-module} \cite{Fai1976}, \textit{Bass module}
\cite{Ann2002}) iff every proper $R$-submodule of $M$ is contained in a
maximal $R$-submodule of $M$, equivalently iff $\mathrm{Rad}(M/L)\neq M/L$
for every $L\lvertneqq _{R}M$;

\emph{f.i.-coatomic} iff $\mathcal{M}^{\mathrm{f.i.}}(L)\neq \varnothing $
for every $L\lvertneqq _{R}^{\mathrm{f.i.}}M$, equivalently iff $_{R}M_{S}$
is coatomic.
\end{definition}

\begin{punto}
By $\mathcal{S}(M)$ ($\mathcal{S}_{\mathrm{f.i.}}(M)$) we denote the \textit{%
possibly empty} class of simple $R$-submodules of $M$ (simple $(R,S)$%
-subbimodules of $_{R}M_{S}$). Let $L\leq _{R}M.$ We set%
\begin{equation*}
\mathcal{S}(L):=\{K\in \mathcal{S}(M)\mid K\subseteq L\}\text{ and }\widehat{%
\mathcal{S}}(L):=\{K\in \mathcal{S}_{\mathrm{f.i.}}(M)\mid K\subseteq L\}.
\end{equation*}%
We say that $L$ is \emph{essential} or \emph{large} in $M,$ and write $%
L\trianglelefteq M$ iff $L\cap \widetilde{L}\neq 0$ for every $0\neq
\widetilde{L}\leq _{R}M.$ The \emph{socle} of $_{R}M$ is defined as
\begin{equation*}
\mathrm{Soc}(M):=\dsum\limits_{L\in \mathcal{S}(M)}L=\dbigcap\limits_{L%
\trianglelefteq M}L\text{ }\;\text{(}:=0\text{ iff }\mathcal{S}%
(M)=\varnothing \text{)}
\end{equation*}
\end{punto}

\begin{definition}
We call $_{R}M:$

\emph{homogenous semisimple} iff $M$ is a (direct) sum of isomorphic simple $%
R$-submodules;

\emph{completely inhomogenous semisimple }iff $M$ is a (direct) sum of
pairwise non-isomorphic simple submodules.
\end{definition}

\begin{definition}
We say $_{R}M$ is

\emph{colocal} (or \textit{cocyclic} \cite{Wis1991}, \textit{subdirectly
irreducible} \cite{AF1974}) iff $M$ contains a \emph{smallest} non-zero $R$%
-submodule that is contained in every non-zero $R$-submodule of $M$,
equivalently iff $\dbigcap\limits_{0\neq L\leq _{R}M}L\neq 0$;

\emph{uniform} iff for any $0\neq L_{1},L_{2}\leq _{R}M,$ also $L_{1}\cap
L_{2}\neq 0$, equivalently iff every non-zero $R$-submodule of $M$ is
essential;

\emph{atomic} iff every $0\neq L\leq _{R}M$ contains a simple $R$-submodule,
equivalently iff $\mathrm{Soc}(L)\neq 0$ for every $0\neq L\leq _{R}M$;

\emph{f.i.-atomic} iff $\widehat{\mathcal{S}}(L)\neq \varnothing $ for every
$0\neq L\leq _{R}^{\mathrm{f.i.}}M$, equivalently iff $_{R}M_{S}$ is atomic.
\end{definition}

\begin{ex}
\label{PG}(\cite[17.13]{Wis1991})\ Let $p\in \mathbb{P}$ and consider the
\emph{Pr\"{u}fer }$p$\emph{-group}%
\begin{equation*}
\mathbb{Z}_{p^{\infty }}:=\sum_{n\in \mathbb{N}}\mathbb{Z}(\frac{1}{p^{n}}+%
\mathbb{Z)=}\bigcup_{n\in \mathbb{N}}\mathbb{Z}(\frac{1}{p^{n}}+\mathbb{Z)}%
\subseteq \mathbb{Q}/\mathbb{Z}.
\end{equation*}%
Every non-zero proper submodule of $\mathbb{Z}_{p^{\infty }}$ is of the form
$\mathbb{Z}(\frac{1}{p^{n}}+\mathbb{Z}\mathfrak{)}$ for some $n\in \mathbb{N}
$ (whence finite) and is fully invariant (i.e. $\mathbb{Z}_{p^{\infty }}$ is
duo). Clearly, $\mathbb{Z}_{p^{\infty }}$ is Artinian but not Noetherian.
Moreover, $\mathbb{Z}_{p^{\infty }}$ is uniserial whence hollow and uniform
(i.e. \emph{biuniform} after \cite{CLVW2006}). Notice that $\mathbb{Z}%
_{p^{\infty }}$ is not local. On the other hand, $\mathbb{Z}_{p^{\infty }}$
is colocal (being an injective hull of $\mathbb{Z}_{p}$).
\end{ex}

\begin{remark}
Hollow modules are not necessarily coatomic (e.g. $\mathbb{Z}_{p^{\infty }}$
is clearly hollow but not coatomic). This shows that the claim of \cite%
{Gon1998} about hollow modules being coatomic is not correct.
\end{remark}

\begin{exs}
\begin{enumerate}
\item By \cite[Proposition 1.14]{Var1979}, $\mathbb{Z}_{p^{\infty }}$ and $%
\mathbb{Z}_{p^{n}}$ -- where $p$ is some prime and $k\in \mathbb{N}$ -- are
the only hollow Abelian groups (up to isomorphism). Since local modules are
precisely the cyclic hollow modules, it follows that the set of local
Abelian groups (up to isomorphism) is $\{\mathbb{Z}_{p^{k}}\mid p\in \mathbb{%
P}$ and $k\in \mathbb{N}\}.$ Any colocal Abelian group is either injective ($%
\simeq \mathbb{Z}_{p^{\infty }}$ for some $p\in \mathbb{P}$) or finite \cite[%
34.14]{Wis1991}.

\item Examples of uniform modules include, in addition to $\mathbb{Z}%
_{p^{\infty }}$ and $\mathbb{Z}_{p^{k}}$ -- where $p\in \mathbb{P}$ and $%
k\in \mathbb{N}$ -- any additive subgroup of $\mathbb{Q}.\mathfrak{\ }$Any
commutative domain $R$ is obviously uniform when considered as an $R$-module
in the canonical way.

\item Let $_{R}M$ be a non-zero uniserial module (e.g. the Abelian group $%
\mathbb{Z}_{p^{\infty }}$). Then $_{R}M$ is trivially hollow and uniform.
Moreover, $_{R}M$ is local (colocal) if and only if $\mathrm{Rad}(M)\neq M$ (%
$\mathrm{Soc}(_{R}M)\neq 0$).

\item The Abelian group $\mathbb{Z}$ is uniform but not atomic. On the other
hand, $\mathbb{Z}$ is coatomic but not hollow. If $M=S_{1}\oplus S_{2},$
where $_{R}S_{1}$ and $_{R}S_{2}$ are simple, then $_{R}M$ is coatomic but
not hollow.
\end{enumerate}
\end{exs}

\begin{definition}
$R$ is said to be

\begin{enumerate}
\item \emph{left Lewy ring }(or \emph{left semiartinian, or left socular }%
\cite[22.10]{Fai1976}) iff every non-zero left $R$-module has a simple
submodule;

\item \emph{left Bass ring }\cite[2.19]{CLVW2006}\emph{\ }(or a \emph{B-ring
}\cite[22.7]{Fai1976})\emph{\ }iff every non-zero left $R$-module has a
maximal submodule.

Right Lewy rings and right Bass rings can be defined analogously. A left and
right Lewy (Bass) rings is said to be a Lewy (Bass) ring.
\end{enumerate}
\end{definition}

\begin{remarks}
\label{LB}For any ring $R,$ one can deduce immediately from the definitions:

\begin{enumerate}
\item $R$ is a left Lewy ring if and only if every non-zero left $R$-module
has essential socle if and only if every non-zero left $R$-module is atomic.

\item $R$ is a left Bass ring if and only if every non-zero left $R$-module
has small radical if and only if every non-zero left $R$-module is coatomic.
\end{enumerate}
\end{remarks}

\begin{ex}
Every right (left) perfect ring is a left Lewy ring \cite{Bac1995} (left
Bass ring \cite{CLVW2006}). Any semiprimary ring is a Bass ring.
\end{ex}

\begin{punto}
$_{R}M$ is said to be \emph{multiplication} iff every $L\leq _{R}M$ is of
the form $L=IM$ for some $I\in \mathcal{I}(R),$ equivalently $L=(L:_{R}M)M$
(i.e. $\mathcal{L}(M)=\mathcal{L}_{m}(M)$). It is obvious that every
multiplication module is duo. Multiplication modules over commutative rings
have been studied intensively in the literature (e.g. \cite{AS2004, PC1995,
Smi1994}). Several results in these papers were generalized to
multiplication modules over rings \textit{close to be commutative (e.g.}
\cite{Tug2003}, \cite{Tug2004}).
\end{punto}

\qquad Recall that a commutative ring $R$ is arithmetical iff $R_{\mathfrak{m%
}}$ is a chain ring for every $\mathfrak{m}\in \mathrm{Max}(R)$ \cite%
{Fuc1949}.

\begin{ex}
A commutative ring $R$ is arithmetical if and only if every finitely
generated ideal of $R$ is multiplication. Let $R$ be arithmetical and $I\in
\mathcal{I}(R)$ be finitely generated. For every $\mathfrak{m}\in \mathrm{Max%
}(R),$ the $R_{\mathfrak{m}}$-ideal $I_{\mathfrak{m}}$ is finitely generated
whence principal since $R_{\mathfrak{m}}$ is a chain ring. Moreover, $R$ is
an fqp-ring, i.e. $_{R}I$ is self-projective, equivalently $_{\overline{R}}I$
is projective \cite{AJK}. Since $I$ is locally principle and $_{\overline{R}%
}I$ is projective, it follows by \cite[Theorem A]{Smi1994} that $_{R}I$ is
multiplication. On the other hand, if every finitely generated ideal of $R$
is multiplication, then every finitely generated ideal of $R$ is locally
principal by \cite[Theorem A]{Smi1994}, whence $R$ is arithmetical by \cite%
{Jen1966}.
\end{ex}

\begin{punto}
$_{R}M$ is said to be \emph{comultiplication} iff every $L\leq _{R}M$ is of
the form $L=(0:_{M}I)$ for some $I\in \mathcal{I}(R),$ equivalently $%
L=(0:_{M}(0:_{R}L))$ (i.e. $\mathcal{L}(M)=\mathcal{L}_{c}(M)$). A ring $R$
for which $_{R}R$ ($R_{R}$) is a comultiplication module is called a \emph{%
left dual} (\emph{right dual}) \textit{ring}. A left dual and right dual
ring is said to be a \emph{dual ring}. For more information on
comultiplication modules and dual rings, the interested reader is referred
to \cite{A-TF2007}, \cite{AS} and \cite{NY2003}.
\end{punto}

\begin{exs}
Let $R$ be left quasi-duo and assume that $PI=IP$ for every $P\in \mathrm{Max%
}(R)$ and $I\in \mathcal{I}(R).$

\begin{enumerate}
\item All multiplication left $R$-modules are coatomic. The proof is similar
to that for multiplication modules over commutative rings \cite[Theorem 2.5
(i)]{A-ES1988} (see also \cite[Proposition 2.4]{Zha2006} for the details in
the non-commutative case).

\item All comultiplication left $R$-modules are atomic. The proof is similar
to that for comultiplication modules over commutative rings \cite[Theorem
3.2 (a)]{A-TF2008}.
\end{enumerate}
\end{exs}

\begin{notation}
For any $L\in \mathrm{Max}(M)$ wet set
\begin{equation*}
L^{e}:=\dbigcap\limits_{K\in \mathrm{Max}(M)\backslash \{L\}}K\text{ \ \ (}%
:=M\text{ iff }\mathrm{Max}(M)=\{L\}\text{).}
\end{equation*}%
Dually, for every $L\in \mathcal{S}(M)$ we set
\begin{equation*}
L_{e}:=\dsum\limits_{K\in \mathcal{S}(M)\backslash \{L\}}K\text{ \ \ (}:=0%
\text{ iff }\mathcal{S}(M)=\{L\}\text{).}
\end{equation*}
\end{notation}

\qquad In \cite{Abu} and \cite{Abu2011}, we introduced the class of modules
with the (\emph{complete}) \emph{max-property} and the class of modules with
the \emph{min-property}. For a study and survey on these modules see \cite%
{Smi-1}:

\begin{punto}
\label{max-prop} We say that $_{R}M$ has the \emph{complete max-property},
iff for any $L\in \mathrm{Max}(M)$ we have $L^{e}\nsubseteqq L.$ We also say
that $_{R}M$ has the \emph{max-property}, iff for any $L\in \mathrm{Max}(M)$
and any finite subset $\mathcal{A}\subseteq \mathrm{Max}(M)\setminus \{L\}$
we have $\bigcap_{K\in \mathcal{A}}K\nsubseteqq L.$
\end{punto}

\begin{lemma}
\label{semi-local}\emph{(\cite[Lemma 3.15]{Abu})} Let $_{R}M$ be
self-projective and duo. Then $M$ has the max-property.
\end{lemma}

\begin{punto}
We say that $_{R}M$ has the \emph{min-property} iff for any simple $R$%
-submodule $L\in \mathcal{S}(M)$ we have $L\nsubseteqq L_{e}.$ Since simple
modules are cyclic, $_{R}M$ has the min-property if and only if for any $%
L\in \mathcal{S}(M)$ and any \textit{finite} subset $\{L_{1},\cdots
,L_{n}\}\subseteq \mathcal{S}(M)\setminus \{L\}$, we have $L\nsubseteqq
\sum_{i=1}^{n}L_{i}.$
\end{punto}

\begin{lemma}
\label{min}\emph{(\cite[Lemma 3.17]{Abu2011})}\ If $_{R}M$ is self-injective
and duo, then $_{R}M$ has the min-property.
\end{lemma}

\subsection*{Topological Spaces}

\qquad In what follows, we fix some definitions and notions for topological
spaces. For further information, the reader might consult any book in
General Topology (e.g. \cite{Bou1966}).

\begin{definition}
We call a topological space $\mathbf{X}$ (\textit{countably}) \textit{compact%
} iff every open cover of $\mathbf{X}$ has a (countable) finite subcover.
Countably compact spaces are also called \textit{Lindel\"{o}f spaces}. Note
that some authors (e.g. \cite{Bou1966, Bou1998}) assume that compact spaces
are in addition Hausdorff.
\end{definition}

\begin{punto}
We say a topological space ${\mathbf{X}}$ is \textit{Noetherian} (\textit{%
Artinian}) iff every ascending (descending) chain of open sets is
stationary, equivalently iff every descending (ascending) chain of closed
sets is stationary.
\end{punto}

\begin{definition}
(e.g. \cite{Bou1966}, \cite{Bou1998}) A non-empty topological space $\mathbf{%
X}$ is called

\begin{enumerate}
\item \emph{ultraconnected} iff the intersection of any two non-empty closed
subsets is non-empty.

\item \emph{irreducible} (or \emph{hyperconnected}) iff $\mathbf{X}$ is not
the union of two proper \textit{closed} subsets, equivalently iff the
intersection of any two non-empty open subsets is non-empty.

\item \emph{connected} iff $\mathbf{X}$ is not the \textit{disjoint} union
of two proper \textit{closed} subsets, equivalently iff the only subsets of $%
\mathbf{X}$ that are \textit{clopen} (i.e. closed and open) are $\varnothing
$ and $\mathbf{X}$.
\end{enumerate}
\end{definition}

\begin{punto}
(\cite{Bou1966}, \cite{Bou1998}) Let $\mathbf{X}$ be a non-empty topological
space. A non-empty subset ${\mathcal{A}}\subseteq {\mathbf{X}}$ is said to
be \textit{irreducible} iff it's an irreducible space w.r.t. the relative
(subspace) topology. A maximal irreducible subspace of $\mathbf{X}$ is
called an \emph{irreducible component}. An irreducible component of a
topological space is necessarily closed. Every irreducible subset of $%
\mathbf{X}$ is contained in an irreducible component of $\mathbf{X},$ whence
$\mathbf{X}$ is the union of its irreducible components. The irreducible
components of a Hausdorff space are just the singleton sets.
\end{punto}

\begin{lemma}
\label{irred}The following are equivalent for ${\mathcal{A}}\subseteq {%
\mathbf{X}}:$

\begin{enumerate}
\item ${\mathcal{A}}$ is irreducible;

\item For any closed subsets $\mathcal{A}_{1},\mathcal{A}_{2}$ of $\mathbf{X}%
:$%
\begin{equation}
{\mathcal{A}}\subseteq \mathcal{A}_{1}\cup \mathcal{A}_{2}\Rightarrow {%
\mathcal{A}}\subseteq \mathcal{A}_{1}\text{ or }{\mathcal{A}}\subseteq
\mathcal{A}_{2}.  \label{irr-clos}
\end{equation}

\item For any open subsets $\mathcal{U}_{1},\mathcal{U}_{2}$ of $\mathbf{X}:$%
\begin{equation}
\mathcal{U}_{1}\cap \mathcal{A}\neq \varnothing \neq \mathcal{U}_{2}\cap
\mathcal{A}\Rightarrow (\mathcal{U}_{1}\cap \mathcal{U}_{2})\cap \mathcal{A}%
\neq \varnothing .  \label{irr-open}
\end{equation}
\end{enumerate}
\end{lemma}

\begin{definition}
Let ${\mathbf{X}}$ be a topological space and $\mathcal{Y}\subseteq {\mathbf{%
X}}$ a closed subset. A point $y\in {\mathcal{Y}}$ is said to be a \textit{%
generic point} iff ${\mathcal{Y}}=\overline{\{y\}}$. If every irreducible
closed subset of ${\mathbf{X}}$ has a unique generic point, then we call ${%
\mathbf{X}}$ a \textit{Sober} space.
\end{definition}

\begin{definition}
A collection $\mathcal{G}$ of subsets of a topological space $\mathbf{X}$ is
\emph{locally finite} iff every point of $\mathbf{X}$ has a neighborhood
that intersects only finitely many elements of $\mathcal{G}.$
\end{definition}

\section{Coprime and second submodules}

\qquad As before, $_{R}M$ is a non-zero left module over the associative
ring $R$ and $S:=\mathrm{End}_{R}(M)^{op}.$ In this section, we introduce
and investigate the spectrum $\mathrm{Spec}^{s}(M)$ of second submodules of $%
M$ and the spectrum $\mathrm{Spec}^{c}(M)$ of coprime submodules of\emph{\ }$%
M.$

\begin{definition}
We call $_{R}M$ (\textit{completely}) \textit{coprime} iff for every $I\in
\mathcal{I}(R)$ ($r\in R$) we have $IM=M$ or $IM=0$ ($rM=M$ or $rM=0$).
Moreover, we say that $K\lvertneqq _{R}M$ is (\textit{completely}) \textit{%
coprime in }$M,$ or a (\textit{completely}) \textit{coprime submodule}, iff
for every $I\in \mathcal{I}(R)$ ($r\in R$) we have $IM+K=M$ or $IM\subseteq
K $ ($rM+K=M$ or $rM\subseteq K$). On the other hand, we say that $0\neq
L\leq _{R}M$ is a (\emph{completely}) \emph{second} \emph{submodule} of $M$
iff $_{R}L$ is a (completely) coprime module.
\end{definition}

\begin{punto}
We set%
\begin{eqnarray*}
\mathrm{Spec}^{\mathrm{c}}(M) &:&=\{K\lvertneqq _{R}M\mid K\text{ is coprime
in }M\}; \\
\mathrm{Spec}^{\mathrm{cc}}(M) &:&=\{K\lvertneqq _{R}M\mid K\text{ is
completely coprime in }M\}.
\end{eqnarray*}%
We say that $_{R}M$ is \textit{coprimeless (c-coprimeless)} iff $\mathrm{Spec%
}^{\mathrm{c}}(M)=\varnothing $ ($\mathrm{Spec}^{\mathrm{cc}}(M)=\varnothing
$). On the other hand, we set%
\begin{eqnarray*}
\mathrm{Spec}^{\mathrm{s}}(M) &:&=\{0\neq K\leq _{R}M\mid \text{ }K\text{ is
a second submodule}\}; \\
\mathrm{Spec}^{\mathrm{cs}}(M) &:&=\{0\neq K\leq _{R}M\mid \text{ }K\text{
is a completely second submodule}\}.
\end{eqnarray*}%
We say that $_{R}M$ is \textit{secondless (}c-\textit{secondless)} iff $%
\mathrm{Spec}^{\mathrm{s}}(M)=\varnothing $ ($\mathrm{Spec}^{\mathrm{cs}%
}(M)=\varnothing $).
\end{punto}

\begin{ex}
(\cite{Smi-2})\ A finitely generated non-zero Abelian group $G$ is coprime
if and only if $G$ is divisible or homogenous semisimple: Let $G$ be
coprime. If $pG=G$ for every $p\in \mathbb{P},$ then $G$ is divisible. If $%
pG=0$ for some $p\in \mathbb{P},$ then $G\simeq \bigoplus \mathbb{Z}/p%
\mathbb{Z}$ is homogenous semisimple. The converse is obvious.
\end{ex}

\begin{remarks}
\label{dr}

\begin{enumerate}
\item $_{R}M$ is (completely) coprime if and only if $0$ is (completely)
coprime in $M$ if and only if $_{R}M$ is a (completely) second submodule of
itself.

\item Let $R$ be commutative. Then $_{R}M$ is coprime if and only if $_{R}M$
is completely coprime.
\end{enumerate}
\end{remarks}

\begin{ex}
Let $R$ be a simple ring that is not a division ring (e.g. $\mathbb{M}%
_{n}(D),$ $n\geq 2,$ the matrix ring of $n\times n$-matrices with entries in
some division ring $D$). Then $_{R}R$ ($R_{R}$) is coprime but not
completely coprime.
\end{ex}

\begin{exs}
\begin{enumerate}
\item \label{ex-hom}The Abelian group $\mathbb{Q}$ is coprime.

\item Let $V$ be a vector space over a division ring $D.$ Then $_{D}V$ is
completely coprime, every $W\lvertneqq _{D}V$ is completely coprime in $V$
and every $0\neq W\leq _{D}V$ is completely second in $V.$

\item Every maximal submodule $L\lvertneqq _{R}M$ is coprime in $M,$ i.e. $%
\mathrm{Max}(M)\subseteq \mathrm{Spec}^{\mathrm{c}}(M).$ In particular, if $%
R $ is a \emph{left max ring,} i.e. a ring over which every non-zero left $R$%
-module has a maximal submodule (e.g. a left perfect ring, or a left
V-ring), then $\mathrm{Spec}^{\mathrm{c}}(M)\neq \varnothing .$

\item Every simple submodule $0\neq L\leq _{R}M$ is second in $M,$ i.e. $%
\mathcal{S}(M)\subseteq \mathrm{Spec}^{\mathrm{s}}(M).$ In particular, $%
\mathrm{Spec}^{\mathrm{s}}(M)\neq \varnothing .$
\end{enumerate}
\end{exs}

\begin{ex}
Every \emph{homogenous} semi-simple $R$-module is clearly coprime. Consider
the semisimple abelian group $M:=\mathbb{Z}/2\mathbb{Z}\oplus \mathbb{Z}/3%
\mathbb{Z}.$ Then $_{\mathbb{Z}}M$ is not coprime. Notice that setting $I:=3%
\mathbb{Z}$ we have $0\neq IM=\mathbb{Z}/2\mathbb{Z}\neq M.$ This shows that
the condition homogenous in cannot be removed.
\end{ex}

\begin{proposition}
\label{IM}\emph{(Compare with \cite{Wij2006, WW2009})} The following are
equivalent:

\begin{enumerate}
\item $_{R}M$ is coprime;

\item $\mathrm{ann}_{R}(M)=\mathrm{ann}_{R}(M/L)$ for every $L\lvertneqq
_{R}^{\mathrm{f.i.}}M;$

\item $\mathrm{ann}_{R}(M)=\mathrm{ann}_{R}(M/L)$ for every $L\lvertneqq
_{R}M;$

\item every $L\lvertneqq _{R}M$ is coprime in $M;$

\item every $L\lvertneqq _{R}^{\mathrm{f.i.}}M$ is coprime in $M;$

\item $R/\mathrm{ann}_{R}(M)$ is cogenerated by $M/L$ for every $L\lvertneqq
_{R}^{\mathrm{f.i.}}M;$

\item $R/\mathrm{ann}_{R}(M)$ is cogenerated by $M/L$ for every $L\lvertneqq
_{R}M.$
\end{enumerate}
\end{proposition}

\begin{Beweis}
$(1)\Rightarrow (2)$ Let $L\lvertneqq _{R}^{\mathrm{f.i.}}M$ and $I:=\mathrm{%
ann}_{R}(M/L).$ Suppose that $\mathrm{ann}_{R}(M)\subsetneqq I.$ Since $%
IM\neq 0,$ it follows that $IM=M,$ a contradiction.

$(2)\Rightarrow (3)$ Let $L\lvertneqq _{R}M$ and $I:=\mathrm{ann}_{R}(M/L).$
Then $I\subseteq \mathrm{ann}_{R}(M/IM)=\mathrm{ann}_{R}(M)$ (notice that $%
IM\lvertneqq _{R}^{\mathrm{f.i.}}M$).

$(3)\Rightarrow (4)$ Let $L\lvertneqq _{R}M$ and $I\in \mathcal{I}(R).$ If $%
\widetilde{L}:=IM+L\subsetneqq M,$ then $I\subseteq \mathrm{ann}_{R}(M/%
\widetilde{L})=\mathrm{ann}_{R}(M)$ whence $IM=0\subseteq L.$

$(4)\Rightarrow (5)$ obvious.

$(5)\Rightarrow (1)$ Set $L:=0\leq _{R}^{\mathrm{f.i.}}M.$

The equivalences $(2)\Leftrightarrow (6)$ and $(3)\Leftrightarrow (7)$
follow from immediately from the fact that, over any ring, a module $N$
cogenerates the ring if and only if it is faithful over it.$\blacksquare $
\end{Beweis}

\begin{ex}
Let $_{R}M$ be \emph{non-Hopf kernel} (i.e. $M/K\simeq M$ for every $%
K\lvertneqq _{R}M$ \cite[HM1987]{HM1987}). Then $_{R}M$ is clearly coprime.
In particular, the Pr\"{u}fer group $\mathbb{Z}_{p^{\infty }}$ is coprime as
as $\mathbb{Z}$-module.
\end{ex}

\qquad The following result follows from the definition and Proposition \ref%
{IM}.

\begin{proposition}
\label{IM+K}The following are equivalent for $K\lvertneqq _{R}M:$

\begin{enumerate}
\item $K$ is (completely) coprime in $M;$

\item $L$ is (completely) coprime in $M$ whenever $K\leq _{R}L\lvertneqq
_{R}^{\mathrm{f.i.}}M;$

\item $L$ is (completely) coprime in $M$ whenever $K\leq _{R}L\lvertneqq
_{R}M;$

\item $M/L$ is a (completely) coprime $R$-module whenever $K\leq
_{R}L\lvertneqq _{R}M;$

\item $M/L$ is a (completely) coprime $R$-module whenever $K\leq
_{R}L\lvertneqq _{R}^{\mathrm{f.i.}}M;$

\item $M/K$ is a (completely) coprime $R$-module.
\end{enumerate}
\end{proposition}

\begin{ex}
(cf. \cite[Example 102]{Ann2002}) Let $R$ be a UFD with an infinite
representative set of the prime elements $\mathbb{P}(R)$ (e.g. $R=\mathbb{Z}$%
). Consider $M=\prod\limits_{p\in \mathbb{P}(R)}R/pR$ and $%
K=\bigoplus\limits_{p\in \mathbb{P}}R/pR.$ Let $r\in R$ and let $%
D_{r}=\{p_{1},\cdots ,p_{k}\}\subset \mathbb{P}(R)$ be a representative set
of the prime divisors of $r.$ Clearly, $rM=\prod\limits_{r\in \mathbb{P}%
(R)\backslash D_{r}}R/pR\nsubseteqq K.$ However, $rM+K=M.$ So, $K$ is
coprime in $M.$
\end{ex}

\begin{remarks}
\begin{enumerate}
\item If $\mathrm{ann}_{R}(M)\in \mathrm{Max}(R),$ then $_{R}M$ is coprime:
Let $I\in \mathcal{I}(R).$ If $IM\neq 0,$ then $\mathrm{ann}_{R}(M)+I=R$
whence $IM=M.$ The converse is not true in general. For example, $_{\mathbb{Z%
}}\mathbb{Q}$ is coprime but $\mathrm{ann}_{\mathbb{Z}}(\mathbb{Q})=0$ is
not a maximal ideal.

\item If $_{R}M_{S}$ is simple, then $_{R}M$ is coprime by Proposition \ref%
{IM}.

\item $R$ is a simple ring if and only if $_{R}R$ ($R_{R}$) is coprime. In
particular, a commutative ring $R$ is a field if and only if $_{R}R$ is
coprime.

\item The following are equivalent:

\begin{enumerate}
\item $R$ is a division ring;

\item $_{R}R$ is completely coprime;

\item $R_{R}$ is completely coprime;

\item $_{R}R$ is coprime and $R$ is (left) duo;

\item $R_{R}$ is coprime and $R$ is (right) duo.
\end{enumerate}
\end{enumerate}
\end{remarks}

\begin{proposition}
\label{c-exact}Let $0\neq L\lvertneqq _{R}M.$ If $K\lvertneqq _{R}M$ is
(completely) coprime in $M,$ then $K\cap L$ is (completely) coprime in $L$
or $(K+L)/L$ is (completely) coprime in $M/L.$
\end{proposition}

\begin{Beweis}
Let $0\neq L\lvertneqq _{R}M$ and assume that $K$ is (completely) coprime in
$M.$

\textbf{Case I. }$K+L\subsetneqq M.$ It follows by Proposition \ref{IM+K}
that $K+L\lvertneqq _{R}M$ is (completely) coprime in $M$ whence $\frac{M/L}{%
(K+L)/L}\simeq M/(K+L)$ is a (completely) coprime $R$-module and so $(K+L)/L$
is (completely) coprime in $M/L.$

\textbf{Case II. }$K+L=M.$ Consider $K\cap L\subsetneqq L.$

For any $I\in \mathcal{I}(R)$ ($r\in R$) we have either $IM+K=M$ ($rM+K=M$)
whence%
\begin{eqnarray*}
IL+(K\cap L) &=&(IL+K)\cap L=(I(L+K)+K)\cap L=(IM+K)\cap L=M\cap L=L \\
\text{(}rL+(K\cap L) &=&(rL+K)\cap L=(r(L+K)+K)\cap L=(rM+K)\cap L=M\cap L=L%
\text{)}
\end{eqnarray*}%
or $IM\subseteq K$ ($rM\subseteq K$) whence $IL\subseteq L\cap K$ ($%
rL\subseteq K\cap L$). Consequently, $K\cap L$ is (completely) coprime in $%
L.\blacksquare $
\end{Beweis}

\begin{corollary}
Let $0\neq L\lvertneqq _{R}M.$ If $L$ and $M/L$ are coprimeless
(c-coprimeless), then $_{R}M$ is also coprimeless (c-coprimeless).
\end{corollary}

\begin{lemma}
\label{cop-prime}

\begin{enumerate}
\item If $_{R}M$ is (completely) coprime, then $\mathrm{ann}_{R}(M)$ is
(completely) prime.

\item If $K\lvertneqq _{R}M$ is (completely) coprime in $M,$ then $\mathfrak{%
p}:=(K:_{R}M)$ is (completely) prime.
\end{enumerate}
\end{lemma}

\begin{Beweis}
\begin{enumerate}
\item Assume that $_{R}M$ is (completely) coprime. Let $IJ\subseteq \mathrm{%
ann}_{R}(M)$ ($ab\in \mathrm{ann}_{R}(M)$) and suppose that $J\nsubseteqq
\mathrm{ann}_{R}(M)$ ($b\notin \mathrm{ann}_{R}(M)$) i.e. $JM\neq 0$ ($%
bM\neq 0$). Since $_{R}M$ is (completely) coprime, we conclude that $JM=M$ ($%
bM=M$) whence $IM=I(JM)=(IJ)M=0$ ($aM=a(bM)=(ab)M=0$) i.e. $I\subseteq
\mathrm{ann}_{R}(M)$ ($a\in \mathrm{ann}_{R}(M)$).

\item The result follows from \textquotedblleft 1\textquotedblright\ and
Proposition \ref{IM+K}.$\blacksquare $
\end{enumerate}
\end{Beweis}

\begin{definition}
Let $K\lvertneqq _{R}M$ be (completely) coprime and consider the
(completely) prime ideal $\mathfrak{p}:=(K:_{R}M).$ We say $K$ is a (\emph{%
completely}) $\mathfrak{p}$\emph{-coprime submodule}\textit{\ of }$M.$
\end{definition}

\qquad The following result extends some results in \cite{A-TF2007} to
comultiplication modules over non-commutative rings and improves some other
results.

\begin{proposition}
\label{c-ann}

\begin{enumerate}
\item Let $_{R}M$ be multiplication. Then $_{R}M$ is coprime if and only if $%
_{R}M$ is simple.

\item Let $_{R}M$ be comultiplication. Then $_{R}M$ is coprime if and only
if $\mathrm{ann}_{R}(M)$ is a prime ideal.
\end{enumerate}
\end{proposition}

\begin{Beweis}
\begin{enumerate}
\item Clearly, every simple module is coprime. Conversely, let $_{R}M$ be
multiplication and coprime. If $N\lvertneqq _{R}M,$ then setting $%
I:=(N:_{R}M)$ we have $IM\neq M$ whence $N=IM=0,$ i.e. $_{R}M$ is simple.

\item If $_{R}M$ is coprime, then $\mathrm{ann}_{R}(M)$ is a prime ideal by
Lemma \ref{cop-prime}. Let $_{R}M$ be comultiplication and assume that $%
\mathrm{ann}_{R}(M)\in \mathrm{Spec}(R).$ Let $I\in \mathcal{I}(R)$ be such
that $IM\neq 0$ and let $J:=(0:_{R}IM).$ Since $JI\subseteq \mathrm{ann}%
_{R}(M)$ and $I\nsubseteqq \mathrm{ann}_{R}(M)$ we obtain $J\subseteq
\mathrm{ann}_{R}(M)$ whence $M=(0:_{M}J)=(0:_{M}(0:_{R}IM))=IM.$
Consequently, $_{R}M$ is coprime.$\blacksquare $
\end{enumerate}
\end{Beweis}

\begin{corollary}
\label{mult-comult}Let $_{R}M$ be multiplication and comultiplication. The
following are equivalent:

\begin{enumerate}
\item $_{R}M$ is coprime;

\item $\mathrm{ann}_{R}(M)$ is a prime ideal;

\item $_{R}M$ is simple.
\end{enumerate}

In particular, if $R$ is a prime ring and $_{R}M$ is faithful,
multiplication and comultiplication, then $_{R}M$ is coprime if and only if $%
_{R}M$ is simple.
\end{corollary}

\begin{punto}
Recall that the ring $R$ is said to be \emph{zero-dimensional} iff every
prime ideal of $R$ is maximal. Examples of zero-dimensional rings include
\emph{biregular rings} \cite[3.18 (6, 7)]{Wis1991} and left (right) perfect
rings. For left (right) duo rings, the notion of zero-dimensionality
coincides with that of $\pi $\emph{-regularity} \cite{Hir1978}. A prime ring
(e.g. a commutative integral domain) is said to be \emph{one-dimensional}
iff every non-zero prime ideal is maximal. In particular, commutative
Dedekind domains are one-dimensional.
\end{punto}

\begin{corollary}
\label{com-s}

\begin{enumerate}
\item If $_{R}M$ is multiplication, then $\mathrm{Spec}^{\mathrm{c}}(M)=%
\mathrm{Max}(M).$

\item If $_{R}M$ is comultiplication, then%
\begin{equation}
\mathrm{Spec}^{\mathrm{s}}(M)=\{0\neq L\leq _{R}M\mid (0:_{R}L)\in \mathrm{%
Spec}(R)\}.  \label{s-prime}
\end{equation}%
If, moreover, $R$ is zero-dimensional, then%
\begin{equation}
\mathrm{Spec}^{\mathrm{s}}(M)=\mathcal{S}(M).  \label{s=s}
\end{equation}
\end{enumerate}
\end{corollary}

\begin{Beweis}
In light of Proposition \ref{c-ann} we need to prove only the last part of
the second statement. Notice that $\mathcal{S}(M)\subseteq \mathrm{Spec}^{%
\mathrm{s}}(M).$ Assume that every prime ideal of $R$ is maximal. Let $K\in
\mathrm{Spec}^{\mathrm{s}}(M)$ so that $(0:_{R}K)\in \mathrm{Spec}(R)$ by
Lemma \ref{cop-prime}, whence a maximal ideal by our assumption on $R.$ It
follows that $K=(0:_{M}(0:_{R}K))$ is simple: if $0\neq K_{1}\lvertneqq
_{R}K,$ then $(0:_{R}K)\subsetneqq (0:_{R}K_{1})\subsetneqq R,$ a
contradiction.$\blacksquare $
\end{Beweis}

\begin{corollary}
\label{r-s-c}

\begin{enumerate}
\item If $R$ is a left due ring (e.g. a commutative ring), then%
\begin{equation}
\mathrm{Spec}^{\mathrm{c}}(_{R}R)=\mathrm{Max}(R).  \label{c=max}
\end{equation}

\item If $R$ is a left dual ring, then%
\begin{equation}
\mathrm{Spec}^{\mathrm{s}}(_{R}R)=\{I\in \mathcal{I}_{l}(R)\mid (0:_{R}I)\in
\mathrm{Spec}(R)\}.
\end{equation}%
If moreover $R$ is zero-dimensional, then%
\begin{equation}
\mathrm{Spec}^{\mathrm{s}}(_{R}R)=\mathrm{Min}(_{R}R).  \label{s=min}
\end{equation}
\end{enumerate}
\end{corollary}

\begin{ex}
A ring $R$ is Quasi-Frobenius if and only if $R$ is dual and Artinian.
Examples of Quasi-Frobenius rings include semisimple Artinian rings, the
group algebra $\mathbb{F}[G]$ where $\mathbb{F}$ is a field and $G$ is a
finite group, and $R/aR$ where $R$ is a commutative PID and $0\neq a\notin
U(R)$ (e.g. $\mathbb{Z}/n\mathbb{Z},$ $n\geq 2$).
\end{ex}

\begin{ex}
Let $B\subseteq \mathbb{Q}$ be the subring consisting of all rational
numbers with odd denominators and $M:=\mathbb{Q}/B.$ Consider the \emph{%
idealization} $R:=B\oplus M$ with multiplication%
\begin{equation}
(b_{1},m_{1})(b_{2},m_{2})=(b_{1}b_{2},b_{1}m_{2}+b_{2}m_{1}).
\end{equation}%
Then $R$ is a dual ring which is not Quasi-Frobenius.
\end{ex}

\begin{remark}
For any module $M,$ the so-called \emph{generalized associated prime} ideals
of $_{R}M$ were introduced in \cite{D-AT2000} as the set%
\begin{equation}
\mathbf{Ass}_{R}(M):=\{\mathfrak{p}\in \mathrm{Spec}(R)\mid \mathfrak{p}%
=(0:_{R}L)\text{ for some }L\leq _{R}M\}.  \label{gen-ass}
\end{equation}%
If $_{R}M$ is comultiplication, then one case easily see that there is a 1-1
correspondence%
\begin{equation}
\mathrm{Spec}^{\mathrm{s}}(M)\longleftrightarrow \mathbf{Ass}_{R}(M),\text{ }%
L\mapsto (0:_{R}L)
\end{equation}%
with inverse $\mathfrak{p}\mapsto (0:_{M}\mathfrak{p}).$ In particular, if $%
R $ is a left dual ring, then there is a 1-1 correspondence%
\begin{equation}
\mathrm{Spec}^{\mathrm{s}}(_{R}R)\longleftrightarrow \mathbf{Ass}_{R}(R).
\end{equation}
\end{remark}

\begin{punto}
For every $L\leq _{R}M$ we define%
\begin{equation*}
\mathcal{V}^{\mathrm{s}}(L):=\{K\in \mathrm{Spec}^{\mathrm{s}}(M)\mid
K\subseteq L\}\;\text{and }\mathcal{X}^{\mathrm{s}}(L):=\{K\in \mathrm{Spec}%
^{\mathrm{s}}(M)\mid K\nsubseteqq L\}.
\end{equation*}%
For every $\mathcal{A}\subseteq \mathrm{Spec}^{\mathrm{s}}(M)$ we set%
\begin{equation}
\mathcal{H}(\mathcal{A}):=\sum\limits_{K\in \mathcal{A}}K\text{ \ (}:=0\text{
iff }\mathcal{A}=\varnothing \text{).}
\end{equation}%
In particular, we set%
\begin{equation}
\mathrm{Corad}_{M}^{\mathrm{s}}(L):=\mathcal{H}(\mathcal{V}^{\mathrm{s}%
}(L))=\sum\limits_{K\in \mathcal{V}^{\mathrm{s}}(L)}K\text{ \ (}:=0\text{
iff }\mathcal{V}^{\mathrm{s}}(L)=\varnothing \text{)}
\end{equation}%
and%
\begin{equation}
\mathrm{Corad}^{\mathrm{s}}(M):=\mathrm{Corad}_{M}^{\mathrm{s}%
}(M)=\sum\limits_{K\in \mathrm{Spec}^{\mathrm{s}}(M)}K\text{ \ (}:=0\text{
iff }\mathrm{Spec}^{\mathrm{s}}(M)=\varnothing \text{).}
\end{equation}%
We say $L\leq _{R}M$ is \textrm{$s$}\emph{-coradical in }$M$ iff $\mathrm{%
Corad}_{M}^{\mathrm{s}}(L)=L.$ In particular, we call $_{R}M$ an \textrm{$s$}%
\emph{-coradical module} iff $\mathrm{Corad}^{\mathrm{s}}(M)=M.$
\end{punto}

\begin{remarks}
\label{corad}

\begin{enumerate}
\item For any $L\leq _{R}M$ we have%
\begin{equation*}
\mathrm{Spec}^{\mathrm{s}}(L)=\mathcal{L}(L)\cap \mathrm{Spec}^{\mathrm{s}%
}(M)=\mathcal{V}^{\mathrm{s}}(L).
\end{equation*}

\item For any $L_{1}\leq _{R}L_{2}\leq _{R}M$ we have $\mathrm{Corad}_{M}^{%
\mathrm{s}}(L_{1})\subseteq \mathrm{Corad}_{M}^{\mathrm{s}}(L_{2}).$
Moreover, for any $L\leq _{R}M$ we have%
\begin{equation*}
\mathrm{Corad}_{M}^{\mathrm{s}}(\mathrm{Corad}_{M}^{\mathrm{s}}(L))=\mathrm{%
Corad}_{M}^{\mathrm{s}}(L).
\end{equation*}

\item Notice that $\mathcal{S}(M)\subseteq \mathrm{Spec}^{\mathrm{s}}(M).$
In particular, if $_{R}M$ is atomic, then for every $0\neq L\leq _{R}M$ we
have $\varnothing \neq \mathcal{S}(L)\subseteq \mathrm{Spec}^{\mathrm{s}}(L)=%
\mathcal{V}^{\mathrm{s}}(L)\subseteq \mathrm{Spec}^{\mathrm{s}}(M).$
\end{enumerate}
\end{remarks}

\begin{definition}
Let $0\neq L\leq _{R}M.$ A maximal element of $\mathcal{V}^{\mathrm{s}}(L),$
if any, is said to be \emph{maximal under }$L.$ A \emph{maximal element of }$%
\mathrm{Spec}^{\mathrm{s}}(M)$ is said to be a \textit{maximal second
submodule} of $M.$
\end{definition}

\begin{lemma}
\label{s-max}Let $_{R}M$ be atomic and comultiplication. For every $0\neq
L\leq _{R}M$ there exists $K\in \mathrm{Spec}^{\mathrm{s}}(M)$ which is
maximal under $L.$
\end{lemma}

\begin{Beweis}
Let $0\neq L\leq _{R}M.$ Since $_{R}M$ is atomic, $\varnothing \neq \mathcal{%
S}(L)\subseteq \mathcal{V}^{\mathrm{s}}(L).$ Let%
\begin{equation*}
K_{1}\subseteq K_{2}\subseteq \cdots \subseteq K_{n}\subseteq
K_{n+1}\subseteq \cdots
\end{equation*}%
be an ascending chain in $\mathcal{V}^{\mathrm{s}}(L)$ and set $\widetilde{K}%
:=\bigcup\limits_{i=1}^{\infty }K_{i}.$ Then we have a descending chain of
prime ideals
\begin{equation}
(0:_{R}K_{1})\supseteq (0:_{R}K_{2})\supseteq \cdots \supseteq
(0:_{R}K_{n})\supseteq (0:_{R}K_{n+1})\supseteq \cdots
\end{equation}%
and it follows that $\mathfrak{p}:=(0:_{R}\widetilde{K})=\bigcap%
\limits_{i=1}^{\infty }(0:_{R}K_{i})$ is a prime ideal. Since $_{R}M$ is
comultiplication, $\widetilde{K}\in \mathcal{V}^{\mathrm{s}}(L)$ by
Corollary \ref{com-s}. By Zorn's Lemma, $\mathcal{V}^{\mathrm{s}}(L)$ has a
maximal element.$\blacksquare $
\end{Beweis}

\begin{punto}
For every $L\leq _{R}M$ we define%
\begin{equation*}
\mathcal{V}^{\mathrm{c}}(L):=\{K\in \mathrm{Spec}^{\mathrm{c}}(M)\mid
L\subseteq K\}\text{ and}\;\mathcal{X}^{\mathrm{c}}(L):=\{K\in \mathrm{Spec}%
^{\mathrm{c}}(M)\mid L\nsubseteqq K\}.
\end{equation*}%
For every $\mathcal{A}\subseteq \mathrm{Spec}^{\mathrm{s}}(M),$ we set%
\begin{equation}
\mathcal{J}(\mathcal{A}):=\bigcap\limits_{K\in \mathcal{A}}K\text{ \ (}:=M%
\text{ iff }\mathcal{A}=\varnothing \text{).}
\end{equation}%
and%
\begin{equation}
\mathrm{Rad}_{M}^{\mathrm{c}}(L):=\mathcal{J}(\mathcal{V}^{\mathrm{c}%
}(L))=\bigcap\limits_{K\in \mathcal{V}^{\mathrm{c}}(L)}K\;\;\text{(}:=M\text{
iff }\mathcal{V}^{\mathrm{c}}(L)=\varnothing \text{).}
\end{equation}%
In particular,%
\begin{equation}
\mathrm{Rad}^{\mathrm{c}}(M):=\mathrm{Rad}_{M}^{\mathrm{c}%
}(0)=\bigcap\limits_{K\in \mathrm{Spec}^{\mathrm{c}}(M)}K\;\;\text{(}:=M%
\text{ iff }\mathrm{Spec}^{\mathrm{c}}(M)=\varnothing \text{).}
\end{equation}%
We say $L\leq _{R}M$ is $\mathrm{c}$\textit{-radical} iff $\mathrm{Rad}_{M}^{%
\mathrm{c}}(L)=L.$
\end{punto}

\begin{remarks}
\label{rad}

\begin{enumerate}
\item For any $L\leq _{R}M$ let%
\begin{equation}
\mathcal{U}(L):=\{\widetilde{L}\leq _{R}M\mid L\leq _{R}\widetilde{L}\leq
_{R}M\}.
\end{equation}%
We have a bijection%
\begin{equation*}
\mathrm{Spec}^{\mathrm{c}}(M)\cap \mathcal{U}(L)\leftrightarrow \mathrm{Spec}%
^{\mathrm{c}}(M/L).
\end{equation*}

\item For any $L_{1}\leq _{R}L_{2}\leq _{R}M$ we have $\mathrm{Rad}_{M}^{%
\mathrm{c}}(L_{1})\subseteq \mathrm{Rad}_{M}^{\mathrm{c}}(L_{2}).$ Moreover,
for any $L\leq _{R}M$ we have%
\begin{equation*}
\mathrm{Rad}_{M}^{\mathrm{c}}(\mathrm{Rad}_{M}^{\mathrm{c}}(L))=\mathrm{Rad}%
_{M}^{\mathrm{c}}(L).
\end{equation*}

\item $\mathrm{Max}(M)\subseteq \mathrm{Spec}^{\mathrm{c}}(M).$ If $_{R}M$
is coatomic, then for every $L\leq _{R}M$ we have $\varnothing \neq \mathcal{%
M}(L)\subseteq \mathcal{V}^{\mathrm{c}}(L)\subseteq \mathrm{Spec}^{\mathrm{c}%
}(M).$
\end{enumerate}
\end{remarks}

\begin{definition}
We call $_{R}M$ (\textit{completely}) \textit{endo-coprime} iff $M_{S}$ is
(completely) coprime. Moreover, we say that $K\lvertneqq _{R}^{\mathrm{f.i.}%
}M$ is a (\textit{completely}) \textit{endo-coprime} $R$\textit{-submodule}
iff $K\leq _{S}M$ is a (completely) coprime submodule.
\end{definition}

\begin{proposition}
\label{endo-c-S}

\begin{enumerate}
\item If $_{R}M$ is (completely)\ endo-coprime, then $S$ is a prime ring (a
domain).

\item Let $M_{S}$ be duo and $_{R}M$ be a self-cogenerator. Then $_{R}M$ is
endo-coprime if and only if $S$ is a prime ring.
\end{enumerate}
\end{proposition}

\begin{Beweis}
\begin{enumerate}
\item If $_{R}M$ is (completely) endo-coprime, then -- by definition -- $%
M_{S}$ is (completely) coprime and it follows by Lemma \ref{cop-prime} that $%
0=\mathrm{ann}_{S}(M)$ is a (completely) prime ideal, i.e. $S$ is a prime
ring (a domain).

\item Assume that $M_{S}$ is duo and $_{R}M$ is a self-cogenerator. Let $%
L\leq _{S}M$ be an arbitrary submodule. Since $M_{S}$ is duo, $L\leq _{B}M$
where $B:=\mathrm{End}(M_{S}).$ Considering the canonical ring morphism $%
\beta :R\rightarrow B,$ we know that $M$ is a $(B,S)$-bimodule and conclude
that $L\leq _{R}^{\mathrm{f.i.}}M$ whence $(0:_{S}L)\in \mathcal{I}(S).$
Since $_{R}M$ is a self-cogenerator, $L=(0:_{M}(0:_{S}L))$ \cite[28.1., 28.2.%
]{Wis1991}. Consequently, $M_{S}$ is a comultiplication module and it
follows by Lemma \ref{c-ann} that $M_{S}$ is coprime, i.e. $_{R}M$ is
endo-coprime.$\blacksquare $
\end{enumerate}
\end{Beweis}

\begin{punto}
We say $_{R}M$ is \textit{divisible} iff $rM=M$ for every $r\in R\backslash
Z(R).$ The sum of all divisible submodules of $_{R}M$ is a divisible
submodule, denoted by $\func{div}(M).$ If $\func{div}(M)=0,$ then $_{R}M$ is
said to be \textit{reduced}. Moreover, $L\leq _{R}M$ is said to be \textit{%
relatively divisible} iff%
\begin{equation*}
rL=rM\cap L\text{ for every }r\in R.
\end{equation*}
\end{punto}

\begin{proposition}
\label{div}Let $K\lvertneqq _{R}L\leq _{R}M.$

\begin{enumerate}
\item Let $_{R}M$ be flat.

\begin{enumerate}
\item If $L\leq _{R}M$ is pure and $K$ is coprime in $M,$ then $K$ is
coprime in $L.$

\item $_{R}M$ is coprime if and only if every non-zero pure submodule of $M$
is second in $M.$
\end{enumerate}

\begin{enumerate}
\item If $L\leq _{R}M$ is relatively divisible and $K$ is completely coprime
in $M,$ then $K$ is completely coprime in $L.$

\item $_{R}M$ is completely coprime if and only if every non-zero relatively
divisible $R$-submodule of $M$ is completely second in $M.$
\end{enumerate}
\end{enumerate}
\end{proposition}

\begin{Beweis}
\begin{enumerate}
\item Let $_{R}M$ be flat.

\begin{enumerate}
\item Assume that $K$ is coprime in $M.$ Since $_{R}M$ is flat and $L\leq
_{R}M$ is pure, we have $IL=IM\cap L$ for every $I\in \mathcal{I}(R)$ (e.g.
\cite[36.6]{Wis1991}). If $IL\nsubseteqq K$ for some $I\in \mathcal{I}(R),$
then indeed $IM\nsubseteqq K$ and so%
\begin{equation}
IL+K=(IM\cap L)+K=(IM+K)\cap L=M\cap L=L.
\end{equation}

\item The proof is similar to that of (a).
\end{enumerate}

\begin{enumerate}
\item Assume that $K$ is coprime in $M.$ Since $L\leq M$ is relatively
divisible, for every $r\in R$ we have
\begin{equation}
rL=rM\cap L\subseteq K\cap L=K
\end{equation}%
or%
\begin{equation*}
rL+K=(rM\cap L)+K=(rM+K)\cap L=M\cap L=L.
\end{equation*}

\item The proof is similar to that of (a).$\blacksquare $
\end{enumerate}
\end{enumerate}
\end{Beweis}

\begin{lemma}
\label{coprime-div}Let $K\lvertneqq _{R}M$ and $\mathfrak{q}:=(K:_{R}M).$
Then $K$ is completely coprime in $M$ if and only if $\mathfrak{q}$ is
completely prime and $M/K$ is a divisible $R/\mathfrak{q}$-module. In
particular, $_{R}M$ is completely coprime if and only if $\overline{R}:=R/%
\mathrm{ann}_{R}(M)$ is a domain and $_{\overline{R}}M$ is divisible.
\end{lemma}

\begin{Beweis}
$(\Rightarrow )$ Assume that $K$ is completely coprime in $M,$ so that $%
\mathfrak{q}$ is completely prime by Lemma \ref{cop-prime} and $\overline{R}%
:=R/\mathfrak{q}$ is a domain. Let $m+K\in M/K$ and consider any $0\neq
\overline{r}\in R/\mathfrak{q}.$ Since $r\notin \mathfrak{q},$ we have $%
rM+K=M$ and so $\overline{r}(M/K)=M/K.$ Consequently, $M/K$ is a divisible $%
R/\mathfrak{q}$-module.

$(\Leftarrow )$ Assume that $\mathfrak{q}$ is completely prime and that $M/K$
is a divisible $R/\mathfrak{q}$-module. Let $r\in R$ be such that $%
rM\nsubseteqq K,$ i.e. $r\notin \mathfrak{q}.$ Since $M/K$ is a divisible $R/%
\mathfrak{q}$-module, we conclude that $\overline{r}(M/K)=M/K$ and so $%
rM+K=M.$ Consequently, $K$ is completely coprime in $M.\blacksquare $
\end{Beweis}

\begin{corollary}
If $K\lvertneqq _{R}M$ and $(K:_{R}M)\in \mathrm{Max}(_{R}R),$ then $K$ is
completely coprime in $M.$ In particular, if $\mathfrak{m}\in \mathcal{I}%
(R)\cap \mathrm{Max}(_{R}R)$ and $\mathfrak{m}M\neq M,$ then $\mathfrak{m}M$
is completely coprime in $M.$
\end{corollary}

\begin{punto}
(\cite{MR1987}, \cite{Mar1972}) Recall that the ring $R$ is said to be
\textit{right} (\textit{left}) \textit{bounded} iff every essential right
(left) ideal of $R$ contains a non-zero \textit{two-sided} ideal. We say
that $R$ is bounded iff $R$ is left and right bounded. Let $R$ be a
Noetherian prime ring with simple Artinian classical ring of quotients $Q.$
If each ideal of $R$ distinct from zero is invertible in $Q,$ then $R$ is
called a \emph{Dedekind prime ring}. For example, all commutative Dedekind
domains and full matrix rings over them are bounded Dedekind prime rings.
Moreover, all (hereditary Noetherian) prime principal ideal rings are
(bounded) Dedekind prime rings.
\end{punto}

\begin{lemma}
\label{bDp}\emph{(\cite[Theorem 3.19]{Mar1972})} Let $R$ be a bounded
Dedekind prime ring. Any $R$-module $M$ possesses a unique largest divisible
submodule $D$ such that $M=D\oplus K$ where $K$ has no divisible submodules.
\end{lemma}

\begin{definition}
Let $R$ be a bounded Dedekind prime ring, $M$ an $R$-module and consider the
decomposition $M=D\oplus K$ from the previous lemma. The submodule $\func{div%
}(M):=D$ ($\mathrm{Red}(M):=K$) is called the \textit{divisible} (\textit{%
reduced}) \textit{part} of $M.$
\end{definition}

\begin{proposition}
\label{div+red}\emph{\ }Let $R$ be a bounded Dedekind prime ring.

\begin{enumerate}
\item If $_{R}M$ is completely coprime, then $_{R}M$ is divisible or reduced.

\item Let $R$ be a domain. Then $\mathrm{Red}(M)$ is completely coprime in $%
M $ if and only if $_{R}M$ is not reduced.
\end{enumerate}
\end{proposition}

\begin{Beweis}
\begin{enumerate}
\item Let $_{R}M$ be completely coprime. If $_{R}M$ is divisible, then we
are done. Suppose that $rM\neq M$ for some $r\in R\backslash Z(R).$ Since $%
_{R}M$ is completely coprime, we have%
\begin{equation*}
0=rM=r(\func{div}(M)\oplus \mathrm{Red}(M))=\func{div}(M)\oplus r\mathrm{Red}%
(M).
\end{equation*}%
It follows that $\mathrm{\func{div}}(M)=0,$ i.e. $_{R}M$ is reduced.

\item Let $R$ be domain. If $\mathrm{Red}(M)$ is completely coprime in $M,$
then in particular $\mathrm{Red}(M)\lvertneqq _{R}M$ and so $_{R}M$ is not
reduced (i.e. $_{R}M$ is divisible or mixed). On the other hand, assume that
$\mathrm{Red}(M)\lvertneqq _{R}M.$ Since $R$ is a domain, we conclude that%
\begin{equation*}
(\mathrm{Red}(M):_{R}M)=(0:_{R}\mathrm{\func{div}}(M))=0.
\end{equation*}%
Applying Lemma \ref{coprime-div}, we conclude that $\mathrm{Red}(M)$ is
completely coprime in $M.\blacksquare $
\end{enumerate}
\end{Beweis}

\begin{corollary}
\begin{enumerate}
\item If $_{R}M$ is c-coprimeless, then $\mathfrak{m}M=M$ for every $%
\mathfrak{m}\in \mathcal{I}(R)\cap \mathrm{Max}(_{R}R).$

\item Let $R$ be a bounded Dedekind prime domain. If $_{R}M$ is
c-coprimeless, then $_{R}M$ is reduced.
\end{enumerate}
\end{corollary}

\begin{proposition}
\label{sum-cop}Let $\{M_{\lambda }\}_{\Lambda }$ be a family of non-zero $R$%
-modules.

\begin{enumerate}
\item We have:

\begin{enumerate}
\item If $\prod\limits_{\lambda \in \Lambda }M_{\lambda }$ is (completely)
coprime, then $_{R}M_{\lambda }$ is (completely) coprime for every $\lambda
\in \Lambda .$

\item If $\bigoplus\limits_{\lambda \in \Lambda }M_{\lambda }$ is
(completely) coprime, then $_{R}M_{\lambda }$ is (completely) coprime for
every $\lambda \in \Lambda .$
\end{enumerate}

\item Assume that $\mathrm{ann}_{R}(M_{\lambda })=\mathfrak{p=}\mathrm{ann}%
_{R}(M_{\gamma })$ for every $\lambda ,\gamma \in \Lambda .$

\begin{enumerate}
\item $\prod\limits_{\lambda \in \Lambda }M_{\lambda }$ is (completely) $%
\mathfrak{p}$-coprime if and only if $_{R}M_{\lambda }$ is (completely) $%
\mathfrak{p}$-coprime for every $\lambda \in \Lambda .$

\item $\bigoplus\limits_{\lambda \in \Lambda }M_{\lambda }$ is (completely) $%
\mathfrak{p}$-coprime if and only if $_{R}M_{\lambda }$ is (completely) $%
\mathfrak{p}$-coprime for every $\lambda \in \Lambda .$
\end{enumerate}
\end{enumerate}
\end{proposition}

\begin{Beweis}
Notice that $\prod\limits_{\lambda \in \Lambda }M_{\lambda }\neq 0\neq
\bigoplus\limits_{\lambda \in \Lambda }M_{\lambda }.$

\begin{enumerate}
\item
\begin{enumerate}
\item Assume that $M:=\prod\limits_{\lambda \in \Lambda }M_{\lambda }$ is
(completely) coprime. For any $\lambda \in \Lambda $ and any $I\in \mathcal{I%
}(R)$ ($r\in R$): If $IM_{\lambda }\neq M_{\lambda }$ ($rM_{\lambda }\neq
M_{\lambda }$), then $IM\neq M$ ($rM\neq M$) and so $IM=0$ ($rM=0$). Whence $%
IM_{\lambda }=0$ ($rM_{\lambda }=0$) for every $\lambda \in \Lambda .$

\item The proof is similar to that of (a).
\end{enumerate}

\item Since $\mathrm{ann}_{R}(M_{\lambda })=\mathfrak{p}=\mathrm{ann}%
_{R}(M_{\gamma })$ for every $\lambda ,\gamma \in \Lambda ,$ we have $%
\mathrm{ann}_{R}(\prod\limits_{\lambda \in \Lambda }M_{\lambda })=\mathfrak{p%
}=\mathrm{ann}_{R}(\bigoplus\limits_{\lambda \in \Lambda }M_{\lambda }).$

\begin{enumerate}
\item Let $M:=\prod\limits_{\lambda \in \Lambda }M_{\lambda }.$ Assume that $%
_{R}M_{\lambda }$ is (completely) $\mathfrak{p}$-coprime for every $\lambda
\in \Lambda .$ For any $I\in \mathcal{I}(R)$ ($r\in R$): If $IM\neq M$ ($%
rM\neq M$), then $IM_{\lambda _{0}}\neq M_{\lambda _{0}}$ ($rM_{\lambda
_{0}}\neq M_{\lambda _{0}}$) for some $\lambda _{0}\in \Lambda .$ Since $%
_{R}M_{\lambda }$ is (completely) $\mathfrak{p}$-coprime, $IM_{\lambda
_{0}}=0$ ($rM_{\lambda _{0}}=0$) whence $I\subseteq \mathfrak{p}$ ($r\in
\mathfrak{p}$) i.e. $IM=0$ ($rM=0$). It follows that $_{R}M$ is (completely)
$\mathfrak{p}$-coprime.

\item The proof is similar to that of (a).$\blacksquare $
\end{enumerate}
\end{enumerate}
\end{Beweis}

\begin{proposition}
\label{sum-less}If $\{M_{\lambda }\}_{\Lambda }$ is a family of coprimeless
(c-coprimeless) $R$-modules, then $\prod\limits_{\lambda \in \Lambda
}M_{\lambda }$ and $\bigoplus\limits_{\lambda \in \Lambda }M_{\lambda }$ are
coprimeless (c-coprimeless).
\end{proposition}

\begin{Beweis}
Let $\{M_{\lambda }\}_{\Lambda }$ be a family of coprimeless (c-coprimeless)
$R$-modules and set $M:=\prod\limits_{\lambda \in \Lambda }M_{\lambda }.$
Suppose that $K\in \mathrm{Spec}^{\mathrm{c}}(M)$ ($K\in \mathrm{Spec}^{%
\mathrm{cc}}(M)$) so that, in particular, $\pi _{\lambda _{0}}(K)\subsetneqq
M_{\lambda _{0}}$ for some $\lambda _{0}\in \Lambda .$ For every $I\in
\mathcal{I}(R)$ ($r\in R$) we have $IM\subseteq K$ ($rM\subseteq K$) whence $%
I\pi _{\lambda _{0}}(M)\subseteq \pi _{\lambda _{0}}(K)$ ($r\pi _{\lambda
_{0}}(M)\subseteq \pi _{\lambda _{0}}(K)$)\ or $IM+K=M$ ($rM+K=M$) whence $%
IM_{\lambda _{0}}+\pi _{\lambda _{0}}(K)=M_{\lambda _{0}}$ ($rM_{\lambda
_{0}}+\pi _{\lambda _{0}}(K)=M_{\lambda _{0}}$). It follows that $\pi
_{\lambda _{0}}(K)\in \mathrm{Spec}^{\mathrm{c}}(M_{\lambda _{0}})$ ($\pi
_{\lambda _{0}}(K)\in \mathrm{Spec}^{\mathrm{cc}}(M_{\lambda _{0}})$), a
contradiction.$\blacksquare $
\end{Beweis}

\qquad The ring $R$ is said to be \textit{binoetherian} (or \textit{weakly
Noetherian }\cite[page 74]{Row2008}) iff $R$ satisfies the ACC on $\mathcal{I%
}(R).$

\begin{ex}
\label{acc-Prop-f}(cf. \cite[Proposition 107]{Ann2002}) Let $R$ be a
binoetherian ring. Then $\mathfrak{q}M$ is coprime in $M$ for some prime
ideal $\mathfrak{q}\in \mathrm{Spec}(R).$ In particular, $\mathrm{Spec}^{%
\mathrm{c}}(M)\neq \varnothing .$
\end{ex}

\begin{Beweis}
Since $R$ is binoetherian, for every $J\in \mathcal{I}(R)$ there exist (e.g.
\cite[Theorem 16.24]{Row2008}, \cite[Theorem 9.2]{Row2006}) $\mathfrak{p}%
_{1},...,\mathfrak{p}_{n}\in \mathrm{Spec}(R)$ such that%
\begin{equation}
\mathfrak{q}_{1}\bullet \cdots \bullet \mathfrak{q}_{k}\subseteq J\subseteq
\mathfrak{q}_{1}\cap \cdots \cap \mathfrak{q}_{k}.  \label{P}
\end{equation}%
In particular, there exist $\mathfrak{p}_{1},...,\mathfrak{p}_{n}\in \mathrm{%
Spec}(R)$ such that $\mathfrak{p}_{1}\bullet \cdots \bullet \mathfrak{p}%
_{n}=0,$ whence%
\begin{equation*}
\mathcal{E}:=\{\mathfrak{p}\in \mathrm{Spec}(R)\mid \mathfrak{p}M\neq
M\}\neq \varnothing .
\end{equation*}%
By assumption, $\mathcal{E}$ has a maximal element $\mathfrak{q}.$ Let $I\in
\mathcal{I}(R).$

\textbf{Case I:}\ $I\subseteq \mathfrak{q}.$ In this case, $IM\subseteq
\mathfrak{q}M.$

\textbf{Case II:} $I\nsubseteqq \mathfrak{q},$ so that $\mathfrak{q}%
\varsubsetneqq J:=I+\mathfrak{q}.$

By (\ref{P}), there exist $\mathfrak{q}_{1},...,\mathfrak{q}_{k}\in \mathrm{%
Spec}(R)$ such that
\begin{equation*}
\mathfrak{q}_{1}\bullet \cdots \bullet \mathfrak{q}_{k}\subseteq J\subseteq
\mathfrak{q}_{1}\cap \cdots \cap \mathfrak{q}_{k}.
\end{equation*}

Since $\mathfrak{q}$ is maximal in $\mathcal{E},$ we have $\mathfrak{q}%
_{j}M=M$ for all $j=1,\cdots ,k,$ whence%
\begin{equation*}
IM+\mathfrak{q}M=JM=M.
\end{equation*}%
Consequently, $\mathfrak{q}M\in \mathrm{Spec^{c}}(M).\blacksquare $
\end{Beweis}

\section{Top$^{\mathrm{s}}$-modules}

\qquad As before $M$ is a non-zero left $R$-module. In this section we
topologize the spectrum of second submodules of $_{R}M$ and investigate the
properties of the induced topology. Several proofs in this section are
similar to proofs of results in \cite{Abu2011}, whence omitted.

\bigskip

\begin{notation}
Set
\begin{equation*}
\begin{tabular}{lllllll}
$\xi ^{\mathrm{s}}(M)$ & $:=$ & $\{\mathcal{V}^{\mathrm{s}}(L)\mid L\in
\mathcal{L}(M)\};$ &  & $\xi _{c}^{\mathrm{s}}(M)$ & $:=$ & $\{\mathcal{V}^{%
\mathrm{s}}(L)\mid L\in \mathcal{L}_{c}(M)\};$ \\
$\tau ^{\mathrm{s}}(M)$ & $:=$ & $\{\mathcal{X}^{\mathrm{s}}(L)\mid L\in
\mathcal{L}(M)\};$ &  & $\tau _{c}^{\mathrm{s}}(M)$ & $:=$ & $\{\mathcal{X}^{%
\mathrm{s}}(L)\mid L\in \mathcal{L}_{c}(M)\};$ \\
$\mathbf{Z}^{\mathrm{s}}(M)$ & $:=$ & $(\mathrm{Spec}^{\mathrm{s}}(M),\tau ^{%
\mathrm{s}}(M));$ &  & $\mathbf{Z}_{c}^{\mathrm{s}}(M)$ & $:=$ & $(\mathrm{%
Spec}^{\mathrm{s}}(M),\tau _{c}^{\mathrm{s}}(M)).$%
\end{tabular}%
\end{equation*}
\end{notation}

\begin{lemma}
\label{s-properties}Consider the class of varieties $\xi ^{\mathrm{s}}(M).$

\begin{enumerate}
\item $\mathcal{V}^{\mathrm{s}}(0)=\varnothing $ and $\mathcal{V}^{\mathrm{s}%
}(M)=\mathrm{Spec}^{\mathrm{s}}(M);$

\item $\bigcap\limits_{\lambda \in \Lambda }\mathcal{V}^{\mathrm{s}%
}(L_{\lambda })=\mathcal{V}^{\mathrm{s}}(\bigcap\limits_{\lambda \in \Lambda
}L_{\lambda })$ for any $\{L_{\lambda }\}_{\Lambda }\subseteq \mathcal{L}%
(M); $

\item For any $I,\widetilde{I}\in \mathcal{I}(R),$ we have
\begin{equation}
\mathcal{V}^{\mathrm{s}}((0:_{M}I))\cup \mathcal{V}^{\mathrm{s}}((0:_{M}%
\widetilde{I}))=\mathcal{V}^{\mathrm{s}}((0:_{M}I)+(0:_{M}\widetilde{I}))=%
\mathcal{V}^{\mathrm{s}}((0:_{M}I\cap \widetilde{I}))=\mathcal{V}^{\mathrm{s}%
}((0:_{M}I\widetilde{I})).  \label{L1+L2}
\end{equation}
\end{enumerate}
\end{lemma}

\begin{Beweis}
Statements \textquotedblleft 1\textquotedblright , \textquotedblleft
2\textquotedblright\ and the inclusions
\begin{equation*}
\mathcal{V}^{\mathrm{s}}((0:_{M}I))\cup \mathcal{V}^{\mathrm{s}}((0:_{M}%
\widetilde{I}))\subseteq \mathcal{V}^{\mathrm{s}}((0:_{M}I)+(0:_{M}%
\widetilde{I}))\subseteq \mathcal{V}^{\mathrm{s}}((0:_{M}I\cap \widetilde{I}%
))\subseteq \mathcal{V}^{\mathrm{s}}((0:_{M}I\widetilde{I}))
\end{equation*}%
in (3)\ are clear. Let $K\in \mathcal{V}^{\mathrm{s}}((0:_{M}I\widetilde{I}%
)) $ and suppose that $K\nsubseteqq (0:_{M}\widetilde{I}),$ whence $%
\widetilde{I}K=K.$ It follows that $IK=I(\widetilde{I}K)=(I\widetilde{I}%
)K=0, $ i.e. $K\subseteq (0:_{M}I).$ Consequently, $K\in \mathcal{V}^{%
\mathrm{s}}((0:_{M}I))\cup \mathcal{V}^{\mathrm{s}}((0:_{M}\widetilde{I}%
)).\blacksquare $
\end{Beweis}

\begin{punto}
We call $0\neq L\leq _{R}M$

\emph{strongly hollow }$\emph{in}$ $M$\emph{\ }\textit{iff for any }$%
L_{1},L_{2}\leq _{R}M$ we have%
\begin{equation}
L\subseteq L_{1}+L_{2}\Rightarrow L\subseteq L_{1}\text{ or }L\subseteq
L_{2};  \label{sh}
\end{equation}

\emph{completely hollow}, iff for any collections $\{L_{\lambda }\}_{\Lambda
}$ of $R$-submodules of $M$ we have:%
\begin{equation}
L=\sum L_{\lambda }\Rightarrow L=L_{\lambda }\text{ for some }\lambda \in
\Lambda ;
\end{equation}
\end{punto}

\begin{remark}
Strongly hollow submodules were considered briefly in \cite{RRW2005} under
the name $\vee $\emph{-coprime submodules. }Completely hollow modules were
introduced under the name \emph{completely coirreducible modules} in \cite%
{A-TF2008}; however, the zero submodule was allowed to be completely
coirreducible which does not fit with our scheme.
\end{remark}

\begin{notation}
We set%
\begin{equation}
\mathcal{SH}(M):=\{L\leq _{R}M\mid \text{ }L\text{ is strongly hollow in }%
M\}.  \label{SH}
\end{equation}
\end{notation}

\begin{remarks}
\label{s-strong}

\begin{enumerate}
\item If $_{R}M$ is uniserial, then every submodule of $M$ is strongly
hollow.

\item If $\mathcal{S}(M)\subseteq \mathcal{SH}(M),$ then $_{R}M$ has the
min-property.
\end{enumerate}
\end{remarks}

\begin{ex}
Let $M$ be an $n$-dimensional vector space over a division ring $D.$ If $%
n\geq 2,$ then $M$ has a vector subspace which is hollow but not strongly
hollow: Let $\mathcal{B}=\{v_{1},\cdots ,v_{n}\}$ be a basis for $V$ and
consider $L=D(v_{1}+\cdots +v_{n}).$ Then $L$ is clearly hollow, being $1$%
-dimensional, but not strongly hollow in $M$ since $L\subseteq Dv_{1}+\cdots
+Dv_{n}$ but $L\nsubseteq Dv_{i}$ for any $i=1,\cdots ,n.$ In particular, $%
\{(x,y)\mid y=x\}$ is hollow but not strongly hollow in $\mathbb{R}^{2}.$ If
$M$ is a uniserial non-zero module with $0\neq L\leq _{R}M$ not finitely
generated, then clearly $L$ is strongly hollow but not completely hollow. In
particular, the Abelian group $\mathbb{Z}_{p^{\infty }}$ is strongly hollow
but not completely hollow.
\end{ex}

\qquad In general, $\xi ^{\mathrm{s}}(M)$ is not closed under finite unions.
This motivates

\begin{definition}
We call $_{R}M$ a \textit{top}$^{\mathrm{s}}$\textit{-module} iff $\xi ^{%
\mathrm{s}}(M)$ is closed under finite unions, equivalently iff $\mathbf{Z}^{%
\mathrm{s}}(M):=(\mathrm{Spec}^{\mathrm{s}}(M),\tau ^{\mathrm{s}}(M))$ is a
topological space.
\end{definition}

\begin{remark}
If $_{R}M$ is secondless (i.e. $\mathrm{Spec}^{\mathrm{s}}(M)=\varnothing $%
), then $_{R}M$ is trivially a top$^{\mathrm{s}}$-module.
\end{remark}

\begin{theorem}
\label{s-top}

\begin{enumerate}
\item $\mathbf{Z}_{c}^{\mathrm{s}}(M):=(\mathrm{Spec}^{\mathrm{s}}(M),\tau
_{c}^{\mathrm{s}}(M))$ is a topological space.

\item If $\mathrm{Spec}^{\mathrm{s}}(M)\subseteq \mathcal{SH}(M),$ then $%
_{R}M$ is a top$^{\mathrm{s}}$-module.
\end{enumerate}
\end{theorem}

\begin{Beweis}
\begin{enumerate}
\item This follows directly from Lemma \ref{s-properties}.

\item This follows from the observation that $\mathrm{Spec}^{\mathrm{s}%
}(M)\subseteq \mathcal{SH}(M)$ if and only if $\mathcal{V}^{\mathrm{s}%
}(L_{1})\cup \mathcal{V}^{\mathrm{s}}(L_{2})=\mathcal{V}^{\mathrm{s}%
}(L_{1}+L_{2})$ (equivalently, $\mathcal{X}^{\mathrm{s}}(L_{1}+L_{2})=%
\mathcal{X}^{\mathrm{s}}(L_{1})\cap \mathcal{X}^{\mathrm{s}}(L_{2})$) for
any $L_{1},L_{2}\leq _{R}M.\blacksquare $
\end{enumerate}
\end{Beweis}

\begin{proposition}
\label{com-prop}Let $_{R}M$ be comultiplication.

\begin{enumerate}
\item Every second submodule of $M$ is strongly hollow (i.e. $\mathrm{Spec}^{%
\mathrm{s}}(M)\subseteq \mathcal{SH}(M)$).

\item Every finitely generated second submodule of $M$ is completely hollow.

\item $_{R}M$ is a top$^{\mathrm{s}}$-module.

\item $_{R}M$ has the min-property.
\end{enumerate}
\end{proposition}

\begin{Beweis}
Let $_{R}M$ be comultiplication.

\begin{enumerate}
\item This follows directly from Lemma \ref{s-properties} and the definition
of comultiplication modules.

\item This follow directly from the definitions and \textquotedblleft
1\textquotedblright .

\item This follows from \textquotedblleft 1\textquotedblright\ and Theorem %
\ref{s-top}.

\item This follows from \textquotedblleft 1\textquotedblright , which yields
$\mathcal{S}(M)\subseteq \mathrm{Spec}^{\mathrm{s}}(M)\subseteq \mathcal{SH}%
(M).\blacksquare $
\end{enumerate}
\end{Beweis}

\begin{lemma}
\label{closure}Let $_{R}M$ be a top$^{\mathrm{s}}$-module. The closure of
any subset $\mathcal{A}\subseteq \mathrm{Spec}^{\mathrm{s}}(M)$ is
\begin{equation}
\overline{\mathcal{A}}=\mathcal{V}^{\mathrm{s}}(\mathcal{H}(\mathcal{A})%
\mathcal{)}.  \label{fc-closure}
\end{equation}
\end{lemma}

\begin{remarks}
\label{simple-char}Let $_{R}M$ be a top$^{\mathrm{s}}$-module and consider
the Zariski topology
\begin{equation}
\mathbf{Z}^{\mathrm{s}}(M):=(\mathrm{Spec}^{\mathrm{s}}(M),\tau ^{\mathrm{s}%
}(M)).
\end{equation}

\begin{enumerate}
\item $\mathbf{Z}^{\mathrm{s}}(M)$ is a $T_{0}$ (Kolmogorov) space.

\item Let $\mathrm{Spec}^{\mathrm{s}}(M)\subseteq \mathcal{SH}(M)$ (e.g. $%
_{R}M$ is comultiplication). Then%
\begin{equation}
\mathcal{B}:=\{\mathcal{X}^{\mathrm{s}}(L)\mid L\leq _{R}M\text{ is finitely
generated}\}
\end{equation}%
is a basis of open sets for $\mathbf{Z}^{\mathrm{s}}(M).$

\item If $L\in \mathrm{Spec}^{\mathrm{s}}(M),$ then $\overline{\{L\}}=%
\mathcal{V}^{\mathrm{s}}(L).$ In particular, for any $K\in \mathrm{Spec}^{%
\mathrm{s}}(M):$%
\begin{equation*}
K\in \overline{\{L\}}\Leftrightarrow K\subseteq L.
\end{equation*}

\item ${\mathcal{X}}^{\mathrm{s}}(L)=\emptyset \Rightarrow \mathrm{Soc}%
(M)\subseteq L.$ The converse holds if, for example, $\mathcal{S}(M)=\mathrm{%
Spec}^{\mathrm{s}}(M).$

\item Let $_{R}M$ be atomic. For every $L\leq _{R}M$ we have ${\mathcal{V}}^{%
\mathrm{s}}(L)=\emptyset $ if and only if $L=0.$

\item Let $0\neq L\leq _{R}M.$ The embedding%
\begin{equation*}
\iota :\mathrm{Spec}^{\mathrm{s}}(L)\rightarrow \mathrm{Spec}^{\mathrm{s}}(M)
\end{equation*}%
induces as continuous map%
\begin{equation*}
\mathbf{\iota }:\mathbf{Z}^{\mathrm{s}}(L)\rightarrow \mathbf{Z}^{\mathrm{s}%
}(M).
\end{equation*}%
This follows from the fact that $\mathbf{\iota }^{-1}(\mathcal{V}^{\mathrm{s}%
}(N))=\mathcal{V}^{\mathrm{s}}(N\cap L)$ for every $R$-submodule $N\leq
_{R}M.$
\end{enumerate}
\end{remarks}

\begin{notation}
Set
\begin{equation*}
\mathbf{CL}(\mathbf{Z}^{\mathrm{s}}(M)):=\{\mathcal{A}\subseteq \mathrm{Spec}%
^{\mathrm{s}}(M)\mid \mathcal{A}=\overline{\mathcal{A}}\}\text{ and }%
\mathcal{CR}^{\mathrm{s}}(M):=\{L\leq _{R}M\mid \mathrm{Corad}_{M}^{\mathrm{s%
}}(L)=L\}.
\end{equation*}
\end{notation}

\begin{theorem}
\label{11}Let $_{R}M$ be a top$^{\mathrm{s}}$-module.

\begin{enumerate}
\item We have an order-preserving bijection%
\begin{equation}
\mathcal{CR}^{\mathrm{s}}(M)\longleftrightarrow \mathbf{CL}(\mathbf{Z}^{%
\mathrm{s}}(M)),\text{ }L\mapsto \mathcal{V}^{\mathrm{s}}(L).  \label{bij}
\end{equation}

\item $\mathbf{Z}^{\mathrm{s}}(M)$ is Noetherian if and only if $_{R}M$
satisfies the DCC condition on $\mathcal{CR}^{\mathrm{s}}(M).$

\item $\mathbf{Z}^{\mathrm{s}}(M)$ is Artinian if and only if $_{R}M$
satisfies the ACC condition on $\mathcal{CR}^{\mathrm{s}}(M).$
\end{enumerate}
\end{theorem}

\begin{theorem}
\label{noth-art}Let $_{R}M$ be a top$^{\mathrm{s}}$-module. If $_{R}M$ is
Artinian \emph{(}Noetherian\emph{)}, then $\mathbf{Z}^{\mathrm{s}}(M)$ is
Noetherian \emph{(}Artinian\emph{)}.
\end{theorem}


\begin{proposition}
\label{duo-irr}Let $_{R}M$ be a top$^{\mathrm{s}}$-module and ${\mathcal{A}}%
\subseteq \mathrm{Spec}^{\mathrm{s}}(M).$

\begin{enumerate}
\item If ${\mathcal{A}}\subseteq \mathrm{Spec}^{\mathrm{s}}(M)$ is
irreducible, then ${\mathcal{H}}({\mathcal{A}})$ is a second submodule of $%
M. $

\item Let $\mathrm{Spec}^{\mathrm{s}}(M)\subseteq \mathcal{SH}(M).$ The
following are equivalent:

\begin{enumerate}
\item ${\mathcal{A}}\subseteq \mathrm{Spec}^{\mathrm{s}}(M)$ is irreducible;

\item ${\mathcal{H}}({\mathcal{A}})$ is a second submodule of $M;$

\item $0\neq {\mathcal{H}}({\mathcal{A}})\leq _{R}M$ is strongly hollow.
\end{enumerate}
\end{enumerate}
\end{proposition}

\begin{Beweis}
\begin{enumerate}
\item Assume that ${\mathcal{A}}$ is irreducible. By definition, ${\mathcal{A%
}}\neq \emptyset $ and so ${\mathcal{H}}({\mathcal{A}})\neq 0.$ Let $I\in
\mathcal{I}(R)$ and suppose that $I{\mathcal{H}}({\mathcal{A}})\neq {%
\mathcal{H}}({\mathcal{A}})$ and $I{\mathcal{H}}({\mathcal{A}})\neq 0.$ Set $%
\mathcal{A}_{1}:=\{K\in \mathcal{A}\mid IK=K\}$ and $\mathcal{A}_{2}:=\{K\in
\mathcal{A}\mid IK=0\}.$ Then ${\mathcal{A}}\subseteq \mathcal{V}^{\mathrm{s}%
}({\mathcal{H}}({\mathcal{A}}_{1}))\cup \mathcal{V}^{\mathrm{s}}({\mathcal{H}%
}(\mathcal{A}_{2})).$ Notice that $\mathcal{A}\nsubseteqq \mathcal{V}^{%
\mathrm{s}}({\mathcal{H}}({\mathcal{A}}_{1}))$ (otherwise, $I{\mathcal{H}}({%
\mathcal{A}})=I{\mathcal{H}}({\mathcal{A}}_{1})={\mathcal{H}}({\mathcal{A}}%
_{1})={\mathcal{H}}({\mathcal{A}})$) and $\mathcal{A}\nsubseteqq \mathcal{V}%
^{\mathrm{s}}({\mathcal{H}}({\mathcal{A}}_{2}))$ (otherwise, $I{\mathcal{H}}(%
{\mathcal{A}})=I{\mathcal{H}}({\mathcal{A}}_{2})=0$), a contradiction.
Consequently, $I{\mathcal{H}}({\mathcal{A}})={\mathcal{H}}({\mathcal{A}})$
or $I{\mathcal{H}}({\mathcal{A}})=0,$ i.e. ${\mathcal{H}}({\mathcal{A}})$ is
a second submodule of $M.$

\item Let $\mathrm{Spec}^{\mathrm{s}}(M)\subseteq \mathcal{SH}(M).$

$(a)\Rightarrow (b)$ follows from \textquotedblleft 1\textquotedblright .

$(b)\Rightarrow (c)$ is obvious.

$(c)\Rightarrow (a)$ Since ${\mathcal{H}}({\mathcal{A}})\neq 0,$ we conclude
that ${\mathcal{A}}\neq \varnothing .$ Suppose that ${\mathcal{A}}\subseteq
\mathcal{V}^{\mathrm{s}}(L_{1})\cup \mathcal{V}^{\mathrm{s}}(L_{2})\subseteq
\mathcal{V}^{\mathrm{s}}(L_{1}+L_{2}),$ i.e. ${\mathcal{H}}({\mathcal{A}}%
)\subseteq L_{1}+L_{2}.$ Since ${\mathcal{H}}({\mathcal{A}})$ is strongly
hollow, we conclude that ${\mathcal{H}}({\mathcal{A}})\subseteq L_{1}$
whence $\mathcal{A}\subseteq \mathcal{V}^{\mathrm{s}}(L_{1})$ or ${\mathcal{H%
}}({\mathcal{A}})\subseteq L_{2}$ whence $\mathcal{A}\subseteq \mathcal{V}^{%
\mathrm{s}}(L_{2}).$ Consequently, $\mathcal{A}$ is irreducible.$%
\blacksquare $
\end{enumerate}
\end{Beweis}

\begin{theorem}
\label{corad-s}Let $_{R}M$ be a top$^{s}$-module.

\begin{enumerate}
\item
\begin{enumerate}
\item If $\mathrm{Spec}^{\mathrm{s}}(M)$ is irreducible, then $\mathrm{Corad}%
_{M}^{\mathrm{s}}(M)$ is a second submodule of $M.$

\item If ${\mathcal{S}}(M)$ is irreducible, then $\mathrm{Soc}(M)$ is a
second submodule of $M.$
\end{enumerate}

\item Let $\mathrm{Spec}^{\mathrm{s}}(M)\subseteq \mathcal{SH}(M).$

\begin{enumerate}
\item The following are equivalent:

\begin{enumerate}
\item $\mathrm{Spec}^{\mathrm{s}}(M)$ is irreducible;

\item $\mathrm{Corad}_{M}^{\mathrm{s}}(M)$ is a second submodule of $M;$

\item $0\neq \mathrm{Corad}_{M}^{\mathrm{s}}(M)\leq _{R}M$ is strongly
hollow.
\end{enumerate}

\item The following are equivalent:

\begin{enumerate}
\item ${\mathcal{S}}(M)$ is irreducible;

\item $\mathrm{Soc}(M)$ is a second submodule of $M;$

\item $0\neq \mathrm{Soc}(M)\leq _{R}M$ is strongly hollow.
\end{enumerate}
\end{enumerate}
\end{enumerate}
\end{theorem}

\begin{ex}
Let $_{R}M$ be a top$^{s}$-module. If $\varnothing \neq {\mathcal{A}}%
\subseteq \mathrm{Spec}^{\mathrm{s}}(M)$ is a chain, then ${\mathcal{A}}$ is
irreducible. In particular, if $_{R}M$ is uniserial, then\textrm{\ }$\mathrm{%
Spec}^{\mathrm{s}}(M)$ is irreducible.
\end{ex}

\begin{notation}
Set%
\begin{equation}
\mathrm{Max}(\mathrm{Spec}^{\mathrm{s}}(M)):=\{K\in \mathrm{Spec}^{\mathrm{s}%
}(M)\mid \text{ }K\text{ is a maximal second submodule of }M\}.
\end{equation}
\end{notation}

\begin{proposition}
\label{max-irr}Let $\mathrm{Spec}^{\mathrm{s}}(M)\subseteq \mathcal{SH}(M)$
(whence $_{R}M$ be a top$^{s}$-module).

\begin{enumerate}
\item The bijection \emph{(\ref{bij})} restricts to a bijection:
\begin{equation}
\mathrm{Spec}^{\mathrm{s}}(M)\longleftrightarrow \{{\mathcal{A}}\mid {%
\mathcal{A}}\subseteq \mathrm{Spec}^{\mathrm{s}}(M)\text{ is an irreducible
closed subset}\}  \label{s-irr-closed}
\end{equation}

\item The bijection \emph{(\ref{s-irr-closed})} restricts to a bijection
\begin{equation*}
\mathrm{Max}(\mathrm{Spec}^{\mathrm{s}}(M))\longleftrightarrow \{{\mathcal{A}%
}\mid {\mathcal{A}}\subseteq \mathrm{Spec}^{\mathrm{s}}(M)\text{ is an
irreducible component}\}.
\end{equation*}
\end{enumerate}
\end{proposition}

\begin{Beweis}
Recall the bijection $\mathcal{CR}^{\mathrm{s}}(M)\overset{\mathcal{V}^{%
\mathrm{s}}(-)}{\longrightarrow }\mathbf{CL}(\mathbf{Z}^{\mathrm{s}}(M)).$

\begin{enumerate}
\item Let $K\in \mathrm{Spec}^{\mathrm{s}}(M).$ Then $K={\mathcal{H}}(%
\mathcal{V}^{\mathrm{s}}(K))$ and so the closed set $\mathcal{V}^{\mathrm{s}%
}(K)$ is irreducible by Proposition \ref{duo-irr} \textquotedblleft
2\textquotedblright . On the other hand, let ${\mathcal{A}}\subseteq \mathrm{%
Spec}^{\mathrm{s}}(M)$ be a closed irreducible subset. Notice that ${%
\mathcal{H}}({\mathcal{A}})$ is second in $M$ by Proposition \ref{duo-irr}
\textquotedblleft 2\textquotedblright\ and that ${\mathcal{A}}=\overline{{%
\mathcal{A}}}=\mathcal{V}^{\mathrm{s}}({\mathcal{H}}({\mathcal{A}})).$

\item Let $K$ be maximal in $\mathrm{Spec}^{\mathrm{s}}(M).$ Then $\mathcal{V%
}^{\mathrm{s}}(K)$ is irreducible by \textquotedblleft 1\textquotedblright .
Let $\mathcal{Y}$ be the irreducible component containing $\mathcal{V}^{%
\mathrm{s}}(K).$ Since $\mathcal{Y}$ is closed, $\mathcal{Y=V}^{\mathrm{s}%
}(L)$ for some $L\in \mathrm{Spec}^{\mathrm{s}}(M).$ Since $\mathcal{V}^{%
\mathrm{s}}(K)\subseteq \mathcal{V}^{\mathrm{s}}(L)$ we have $K\subseteq L.$
Since $K\in \mathrm{Max}(\mathrm{Spec}^{\mathrm{s}}(M)),$ we conclude that $%
K=L$ and so $\mathcal{V}^{\mathrm{s}}(K)$ is an irreducible component of $%
\mathrm{Spec}^{\mathrm{s}}(M).$

Conversely, let $\mathcal{Y}$ be an irreducible component of $\mathrm{Spec}^{%
\mathrm{s}}(M).$ Since $\mathcal{Y}$ is closed and irreducible, it follows
by \textquotedblleft 1\textquotedblright\ that $\mathcal{Y}=\mathcal{V}^{%
\mathrm{s}}(L)$ for some $L\in \mathrm{Spec}^{\mathrm{s}}(M).$ Suppose that $%
L$ is not maximal in $\mathrm{Spec}^{\mathrm{s}}(M),$ i.e. there exists $%
K\in \mathrm{Spec}^{\mathrm{s}}(M)$ such that $L\subsetneqq K\subseteq M.$
It follows that $\mathcal{V}^{\mathrm{s}}(L)\subsetneqq \mathcal{V}^{\mathrm{%
s}}(K),$ a contradiction since $\mathcal{V}^{\mathrm{s}}(K)\subseteq \mathrm{%
Spec}^{\mathrm{s}}(M)$ is irreducible by \textquotedblleft
1\textquotedblright . We conclude that $L$ is maximal in $\mathrm{Spec}^{%
\mathrm{s}}(M).\blacksquare $
\end{enumerate}
\end{Beweis}

\begin{corollary}
If $\mathrm{Spec}^{\mathrm{s}}(M)\subseteq \mathcal{SH}(M),$ then $\mathrm{%
Spec}^{\mathrm{s}}(M)$ is a Sober space.
\end{corollary}

\begin{Beweis}
Let ${\mathcal{A}}\subseteq \mathrm{Spec}^{\mathrm{s}}(M)$ be an irreducible
closed subset. By Proposition \ref{max-irr} \textquotedblleft
1\textquotedblright , ${\mathcal{A}}={\mathcal{V}}^{\mathrm{s}}(K)$ for some
$K\in \mathrm{Spec}^{\mathrm{s}}(M)$. It follows that
\begin{equation*}
{\mathcal{A}}=\overline{{\mathcal{A}}}={\mathcal{V}}^{\mathrm{s}}({\mathcal{H%
}}({\mathcal{A}}))={\mathcal{V}}^{\mathrm{s}}(K)=\overline{\{K\}},
\end{equation*}%
i.e. $K$ is a generic point for ${\mathcal{A}}$. If $L$ is a generic point
of ${\mathcal{A}}$, then ${\mathcal{V}}^{\mathrm{s}}(K)={\mathcal{V}}^{%
\mathrm{s}}(L)$ whence $K=L.\blacksquare $
\end{Beweis}

\begin{proposition}
\label{uniform}Let $_{R}M$ be an atomic top$^{\mathrm{s}}$-module. Then $%
_{R}M$ is uniform if and only if $\mathrm{Spec}^{\mathrm{s}}(M)$ is
ultraconnected.
\end{proposition}

\begin{theorem}
\label{compact}Let $_{R}M$ be an atomic top$^{\mathrm{s}}$-module.

\begin{enumerate}
\item If $\mathcal{S}(M)$ is countable, then $\mathbf{Z}^{\mathrm{s}}(M)$ is
countably compact.

\item If $\mathcal{S}(M)$ is finite, then $\mathbf{Z}^{\mathrm{s}}(M)$ is
compact.
\end{enumerate}
\end{theorem}

\begin{ex}
\label{pr-compact}The Pr\"{u}fer group $\mathbb{Z}(p^{\infty })$ is an
atomic top$^{\mathrm{s}}$-module over $\mathbb{Z}.$ Since $\mathcal{S}(%
\mathbb{Z}(p^{\infty })=\{\mathbb{Z}(\frac{1}{p}+\mathbb{Q}/\mathbb{Z})\}$
is finite, we conclude that $\mathbf{Z}^{\mathrm{s}}(\mathbb{Z}(p^{\infty
})) $ is compact.
\end{ex}

\begin{proposition}
\label{it-irr}Let $_{R}M$ be a top$^{\mathrm{s}}$-module and assume that
every second submodule of $M$ is simple.

\begin{enumerate}
\item If $_{R}M$ has the min-property, then $\mathrm{Spec}^{\mathrm{s}}(M)$
is discrete.

\item $M$ has a unique simple $R$-submodule if and only if $_{R}M$ has the
min-property and $\mathrm{Spec}^{\mathrm{s}}(M)$ is connected.
\end{enumerate}
\end{proposition}

\begin{Beweis}
\begin{enumerate}
\item If $_{R}M$ has the min-property, then for every $K\in \mathrm{Spec}^{%
\mathrm{s}}(M)={\mathcal{S}}(M)$ we have $\{K\}={\mathcal{X}}(\{K\}_{e})$ an
open set. Since every singleton set is open, $\mathrm{Spec}^{\mathrm{s}}(M)$
is discrete.

\item ($\Rightarrow $) Assume that $_{R}M$ has a unique simple $R$%
-submodule. Clearly, $_{R}M$ has the min-property and $\mathrm{Spec}^{%
\mathrm{s}}(M)$ is connected since it consists of only one point.

($\Leftarrow $) Assume that $_{R}M$ has the min-property and that $\mathrm{%
Spec}^{\mathrm{s}}(M)$ is connected. By \textquotedblleft
1\textquotedblright , $\mathrm{Spec}^{\mathrm{s}}(M)$ is discrete and so ${%
\mathcal{S}}(M)=\mathrm{Spec}^{\mathrm{s}}(M)$ has only one point since a
discrete connected space cannot contain more than one-point.$\blacksquare $
\end{enumerate}
\end{Beweis}

\begin{theorem}
\label{count-compact}Let $_{R}M$ be an atomic top$^{\mathrm{s}}$-module and
assume that every second submodule of $M$ is simple. If $_{R}M$ has the
min-property, then

\begin{enumerate}
\item $\mathrm{Spec}^{\mathrm{s}}(M)$ is countably compact if and only if ${%
\mathcal{S}}(M)$ is countable.

\item $\mathrm{Spec}^{\mathrm{s}}(M)$ is compact if and only if ${\mathcal{S}%
}(M)$ is finite.
\end{enumerate}
\end{theorem}

\qquad As a direct consequence of Proposition \ref{it-irr} we obtain:

\begin{theorem}
\label{colocal}Let $_{R}M$ be atomic and assume that $\mathcal{S}(M)=\mathrm{%
Spec}^{\mathrm{s}}(M)\subseteq \mathcal{SH}(M)$ so that $M$ is a top$^{%
\mathrm{s}}$-module. Then $_{R}M$ is colocal if and only if $\mathrm{Spec}^{%
\mathrm{s}}(M)$ is connected.
\end{theorem}

\begin{lemma}
\label{1n}Let $\mathrm{Spec}^{\mathrm{s}}(M)\subseteq \mathcal{SH}(M).$ If $%
n\geq 2$ and $\mathcal{A}=\{K_{1},...,K_{n}\}\subseteq \mathrm{Spec}^{%
\mathrm{s}}(M)$ is a connected subset, then for every $i\in \{1,...,n\},$
there exists $j\in \{1,...,n\}\backslash \{i\}$ such that $K_{i}\leq
_{R}K_{j}$ or $K_{j}\leq _{R}K_{i}.$
\end{lemma}

\begin{proposition}
\label{lf}Let $_{R}M$ be an atomic top$^{s}$-module with the min-property
and let $\varnothing \neq \mathcal{A}=\{K_{\lambda }\}_{\Lambda }\subseteq
\mathcal{S}(M).$ If $\left\vert \mathcal{S}(L)\right\vert <\infty $ for
every $L\in \mathrm{Spec}^{\mathrm{s}}(M),$ then $\mathcal{A}$ is locally
finite.
\end{proposition}

\begin{lemma}
\label{s-t1}Let $_{R}M$ be an atomic top$^{\mathrm{s}}$-module. Then the
following are equivalent for any $L\leq _{R}M:$

\begin{enumerate}
\item $L \in {\mathcal{S}}(M)$;

\item $L$ is a second submodule of $M$ and ${\mathcal{V}}^{\mathrm{s}%
}(L)=\{L\}$;

\item $\{L\}$ is closed in ${\mathbf{Z}}_{M}^{\mathrm{s}}$.
\end{enumerate}
\end{lemma}

\begin{proposition}
\label{T1} If $_{R}M$ is an atomic top$^{\mathrm{s}}$-module, then $\mathrm{%
Spec}^{\mathrm{s}}(M)={\mathcal{S}}(M)$ if and only if $\mathbf{Z}^{\mathrm{s%
}}(M)$ is $T_{1}$ \emph{(}Fr\'{e}cht space\emph{)}.
\end{proposition}

\qquad Combining Propositions \ref{it-irr} and \ref{T1} we obtain

\begin{theorem}
\label{T2}Let $_{R}M$ be an atomic top$^{\mathrm{s}}$-module with the
min-property. The following are equivalent:

\begin{enumerate}
\item $\mathrm{Spec}^{\mathrm{s}}(M)=\mathcal{S}(M);$

\item $\mathbf{Z}^{\mathrm{s}}(M)$ is discrete;

\item $\mathbf{Z}^{\mathrm{s}}(M)$ is $T_{2}$ \emph{(}Hausdorff space\emph{)}%
;

\item $\mathbf{Z}^{\mathrm{s}}(M)$ is $T_{1}$ \emph{(}Fr\'{e}cht space\emph{)%
}.
\end{enumerate}
\end{theorem}

\section{Top$^{\mathrm{c}}$-modules}

\qquad As before, $_{R}M$ is a non-zero left $R$-module. In this section we
topologize the spectrum of coprime submodules of $_{R}M$ and investigate the
properties of the induced topology. Several proofs in this section are
similar to proofs of similar results in \cite{Abu} and dual results in
Section 4, whence omitted.

\begin{notation}
Recall that for every $L\leq _{R}M$ we define%
\begin{equation*}
\mathcal{V}^{\mathrm{c}}(L):=\{K\in \mathrm{Spec}^{\mathrm{c}}(M)\mid
L\subseteq K\}\text{ and}\;\mathcal{X}^{\mathrm{c}}(L):=\{K\in \mathrm{Spec}%
^{\mathrm{c}}(M)\mid L\nsubseteqq K\}.
\end{equation*}%
Moreover, we set%
\begin{equation*}
\begin{tabular}{lllllll}
$\xi ^{\mathrm{c}}(M)$ & $:=$ & $\{\mathcal{V}^{\mathrm{c}}(L)\mid L\in
\mathcal{L}(M)\};$ &  & $\xi _{m}^{\mathrm{c}}(M)$ & $:=$ & $\{\mathcal{V}^{%
\mathrm{c}}(L)\mid L\in \mathcal{L}_{m}(M)\};$ \\
$\tau ^{\mathrm{c}}(M)$ & $:=$ & $\{\mathcal{X}^{\mathrm{c}}(L)\mid L\in
\mathcal{L}(M)\};$ &  & $\tau _{m}^{\mathrm{c}}(M)$ & $:=$ & $\{\mathcal{X}^{%
\mathrm{c}}(L)\mid L\in \mathcal{L}_{m}(M)\};$ \\
$\mathbf{Z}^{\mathrm{c}}(M)$ & $:=$ & $(\mathrm{Spec}^{\mathrm{c}}(M),\tau ^{%
\mathrm{c}}(M));$ &  & $\mathbf{Z}_{m}^{\mathrm{c}}(M)$ & $:=$ & $(\mathrm{%
Spec}^{\mathrm{c}}(M),\tau _{m}^{\mathrm{c}}(M)).$%
\end{tabular}%
\end{equation*}
\end{notation}

\begin{lemma}
\label{c-properties}Consider the class of varieties $\xi ^{\mathrm{c}}(M).$

\begin{enumerate}
\item $\mathcal{V}^{\mathrm{c}}(M)=\varnothing $ and $\mathcal{V}^{\mathrm{c}%
}(0)=\mathrm{Spec}^{\mathrm{c}}(M);$

\item $\bigcap\limits_{\lambda \in \Lambda }\mathcal{V}^{\mathrm{c}%
}(L_{\lambda })=\mathcal{V}^{\mathrm{c}}(\sum\limits_{\lambda \in \Lambda
}L_{\lambda })$ for any $\{L_{\lambda }\}_{\Lambda }\subseteq \mathcal{L}%
(M); $

\item For any $I,\widetilde{I}\in \mathcal{I}(R),$ we have
\begin{equation*}
\mathcal{V}^{\mathrm{c}}(IM)\cup \mathcal{V}^{\mathrm{c}}(\widetilde{I}M)=%
\mathcal{V}^{\mathrm{c}}(IM\cap \widetilde{I}M)=\mathcal{V}^{\mathrm{c}%
}((I\cap \widetilde{I})M)=\mathcal{V}^{\mathrm{c}}(I\widetilde{I}M).
\end{equation*}
\end{enumerate}
\end{lemma}

\begin{Beweis}
Statements \textquotedblleft 1\textquotedblright , \textquotedblleft
2\textquotedblright\ and the inclusions
\begin{equation}
\mathcal{V}^{\mathrm{c}}(IM)\cup \mathcal{V}^{\mathrm{c}}(\widetilde{I}%
M)\subseteq \mathcal{V}^{\mathrm{c}}(IM\cap \widetilde{I}M)\subseteq
\mathcal{V}^{\mathrm{c}}((I\cap \widetilde{I})M)\subseteq \mathcal{V}^{%
\mathrm{c}}(I\widetilde{I}M)
\end{equation}%
in (3)\ are clear. Let $K\in \mathcal{V}^{\mathrm{c}}(I\widetilde{I}M)$ and
suppose that $K\notin \mathcal{V}^{\mathrm{c}}(\widetilde{I}M).$ Then $M=%
\widetilde{I}M+K,$ whence $IM=I(\widetilde{I}M+K)=I(\widetilde{I}M)+IK=(I%
\widetilde{I})M+IK\subseteq K,$ i.e. $K\in \mathcal{V}^{\mathrm{c}%
}(IM)).\blacksquare $
\end{Beweis}

In general, $\xi ^{\mathrm{c}}(M)$ is not closed under finite unions. This
motivates

\begin{definition}
We call $_{R}M$ a \textit{top}$^{\mathrm{c}}$\textit{-module} iff $\xi ^{%
\mathrm{c}}(M)$ is closed under finite unions, equivalently $\mathbf{Z}^{%
\mathrm{c}}(M):=(\mathrm{Spec}^{\mathrm{c}}(M),\tau ^{\mathrm{c}}(M))$ is a
topological space.
\end{definition}

\begin{remark}
If $_{R}M$ is coprimeless, i.e. $\mathrm{Spec}^{\mathrm{c}%
}(_{R}M)=\varnothing ,$ then $_{R}M$ is trivially a \textit{top}$^{\mathrm{c}%
}$\textit{-module.}
\end{remark}

\begin{ex}
(cf. \cite[Remark 105]{Ann2002}) Let $R$ have a unique prime ideal $%
\mathfrak{p}.$ For any $_{R}M$ we have%
\begin{equation}
\mathrm{Spec}^{\mathrm{c}}(_{R}M)=\left\{
\begin{array}{ccc}
\mathfrak{p}M, &  & \mathfrak{p}M\neq M \\
&  &  \\
\varnothing , &  & \mathfrak{p}M=M%
\end{array}%
\right.
\end{equation}%
So, if $\mathfrak{p}M=M,$ then $M$ is a \textit{top}$^{\mathrm{c}}$\textit{%
-module. }
\end{ex}

\begin{ex}
(\cite[Example 106]{Ann2002}) Let $R=\mathbb{Q}[x_{1},x_{2},\cdots ],$ where
$x_{i}^{2}=0$ for every $i\in \mathbb{N},$ be the commutative local ring
with \emph{unique prime ideal} $\mathfrak{p}=(x_{1},x_{2},\cdots ).$ Let $%
E=\{e_{1},e_{2},\cdots \}$ be a countably infinite set, $F$ the free $R$%
-module with basis $E$ and $M$ the $R$-module $F$ modulo the relations $%
e_{i}x_{i}=e_{i-1}$ for each $i\geq 2.$ Indeed $\mathfrak{p}M=M,$ whence $%
\mathrm{Spec}^{\mathrm{c}}(M)=\varnothing .$ So, $M$ is a \textit{top}$^{%
\mathrm{c}}$\textit{-module.}
\end{ex}

\begin{ex}
(\cite{Smi-2}) Let $R$ have a unique (completely) prime ideal $\mathfrak{m}.$
If $\mathfrak{m}$ is idempotent, then $\mathfrak{m}$ is (c-coprimeless)
coprimeless: Suppose that $\mathfrak{m}$ contains a left subideal $%
J\varsubsetneqq \mathfrak{m}$ that is (completely) coprime in $\mathfrak{m}.$
By Lemma \ref{coprime-div}, $\mathfrak{p}:=(J:_{R}\mathfrak{m})$ is a
(completely) prime ideal, whence $\mathfrak{p}=\mathfrak{m}$ and
consequently $\mathfrak{m}=\mathfrak{m}^{2}=\mathfrak{pm}\subseteq J,$ a
contradiction. Consequently, $\mathfrak{m}$ is (c-coprimeless) coprimeless.
We conclude that if $\mathfrak{m}$ is idempotent, then $_{R}\mathfrak{m}$ is
a \textit{top}$^{\mathrm{c}}$\textit{-module.}
\end{ex}

\begin{ex}
Let $R:=\mathbb{F}[\mathbb{G}],$ where $\mathbb{F}$ is a field, $G=\mathbb{Z}%
(p^{\infty })$ and $\mathfrak{m}=\dsum\limits_{g\in G}R(g-1):$ Notice that $%
R/\mathfrak{m}\simeq \mathbb{F}$ whence $\mathfrak{m}\in \mathrm{Max}(R).$
For every $x\in G,$ there exists $y\in G$ such that $x=y^{p}.$ So $%
(x-1)=(y-1)^{p}\in I^{p}\subseteq \mathfrak{m}^{2}.$ Therefore, $\mathfrak{m}%
=\mathfrak{m}^{2}$ and so $\mathfrak{m}$ is the only prime ideal of $R.$
Consequently, $_{R}\mathfrak{m}$ is a \textit{top}$^{\mathrm{c}}$\textit{%
-module.}
\end{ex}

\begin{punto}
We call $L\lvertneqq _{R}M:$

\emph{irreducible} iff for any $L_{1},L_{2}\leq _{R}M:$%
\begin{equation}
L_{1}\cap L_{2}=L\Rightarrow L_{1}=L\text{ or }L_{2}=L;  \label{irr}
\end{equation}%
\emph{strongly irreducible}\textit{\ }iff for any $L_{1},L_{2}\leq _{R}M:$%
\begin{equation}
L_{1}\cap L_{2}\subseteq L\Rightarrow L_{1}\subseteq L\text{ or }%
L_{2}\subseteq L.  \label{s-irr}
\end{equation}%
\emph{completely irreducible}, iff any collections $\{L_{\lambda
}\}_{\Lambda }$ of $R$-submodules of $_{R}M$ we have:%
\begin{equation}
\dbigcap\limits_{\lambda \in \Lambda }L_{\lambda }=L\Rightarrow L_{\lambda
}=L\text{ for some }\lambda \in \Lambda ;
\end{equation}
\end{punto}

\begin{remark}
Strongly irreducible ideals were introduced first by Bourbaki \cite[p. 301,
Exercise 34]{Bou1998} and named \emph{quasi-prime ideals}. Recently, they
were investigated by Heinzer et al. \cite{HRR2002} while the class of
completely irreducible ideals was introduced by Fuchs et al. in \cite%
{FHO2006}. The corresponding notions for submodules of modules over
commutative rings were investigated in \cite{E-A2005} and \cite{A-TF2008},
respectively.
\end{remark}

\begin{notation}
We set%
\begin{equation}
\mathcal{SI}(M):=\{L\leq _{R}M\mid \text{ }L\text{ is strongly irreducible}%
\}.  \label{SI(M)}
\end{equation}
\end{notation}

\begin{exs}
\label{strong}

\begin{enumerate}
\item If $_{R}M$ is uniserial, then every submodule of $M$ is strongly
irreducible.

\item $\mathrm{Spec}^{\mathrm{c}}(M)\subseteq \mathcal{SI}(M)$ if and only
if $\mathcal{V}^{\mathrm{c}}(L_{1})\cup \mathcal{V}^{\mathrm{c}}(L_{2})=%
\mathcal{V}^{\mathrm{c}}(L_{1}\cap L_{2})$ (equivalently, $\mathcal{X}^{%
\mathrm{c}}(L_{1})\cap \mathcal{X}^{\mathrm{c}}(L_{2})=\mathcal{X}^{\mathrm{c%
}}(L_{1}\cap L_{2})$) for all $L_{1},L_{2}\leq _{R}M.$

\item If $\mathrm{Max}(M)\subseteq \mathcal{SI}(M),$ then $_{R}M$ has the
max-property.
\end{enumerate}
\end{exs}

\qquad As a direct consequence of Lemma \ref{c-properties} we obtain

\begin{theorem}
\label{c-top}

\begin{enumerate}
\item $\mathbf{Z}_{m}^{\mathrm{c}}(M):=(\mathrm{Spec}^{\mathrm{c}}(M),\tau
_{m}^{\mathrm{c}}(M))$ is a topological space.

\item If $\mathrm{Spec}^{\mathrm{c}}(M)\subseteq \mathcal{SI}(M),$ then $%
_{R}M$ is a top$^{\mathrm{c}}$ module.
\end{enumerate}
\end{theorem}

\qquad Recall that a left $R$-module $N$ is said to be \emph{finitely
cogenerated} iff for any monomorphism $N\overset{f}{\longrightarrow }%
\dprod_{\lambda \in \Lambda }N_{\lambda },$ there exists a finite subset $%
\{\lambda _{1},\cdots ,\lambda _{n}\}\subseteq \Lambda $ such that%
\begin{equation*}
N\overset{f}{\longrightarrow }\dprod_{\lambda \in \Lambda }N_{\lambda }%
\overset{\pi }{\longrightarrow }\dbigoplus\limits_{i=1}^{n}N_{\lambda _{i}}
\end{equation*}%
is injective.

\begin{proposition}
\label{mul-prop}Let $_{R}M$ be a multiplication module.

\begin{enumerate}
\item Every coprime submodule of $M$ is strongly irreducible (i.e. $\mathrm{%
Spec}^{\mathrm{c}}(M)\subseteq \mathcal{SI}(M)$).

\item If $L\in \mathrm{Spec}^{\mathrm{c}}(M)$ is such that $M/L$ is finitely
cogenerated, then $L$ is completely irreducible.

\item $_{R}M$ is a top$^{\mathrm{c}}$ module.

\item $_{R}M$ has the max-property.
\end{enumerate}
\end{proposition}

\begin{Beweis}
\begin{enumerate}
\item This follows directly from Lemma \ref{c-properties}.

\item This follows directly from the definitions and \textquotedblleft
1\textquotedblright .

\item This follows directly from \textquotedblleft 1\textquotedblright\ and
Theorem \ref{c-top}.

\item This follows directly from \textquotedblleft 1\textquotedblright ,
which yields $\mathrm{Max}(M)\subseteq \mathrm{Spec}^{\mathrm{c}%
}(M)\subseteq \mathcal{SI}(M).\blacksquare $
\end{enumerate}
\end{Beweis}

\qquad The following example provides a top$^{\mathrm{c}}$-module which is
not multiplication.

\begin{ex}
(\cite{Smi-2}) Consider the Abelian group $G:=\bigoplus\limits_{p\in \mathbb{%
P}}\mathbb{Z}/p\mathbb{Z}.$ Let $K\lvertneqq _{R}G$ be coprime in $G.$
Notice that $\mathfrak{p}:=(K:G)\neq 0$ (otherwise, $G/K$ would be a
divisible $\mathbb{Z}$-module by Lemma \ref{coprime-div}, a contradiction).
So, $\mathfrak{p}=q\mathbb{Z}$ for some $q\in \mathbb{P}$ and $%
K=\bigoplus\limits_{p\in \mathbb{P}\backslash \{q\}}\mathbb{Z}/p\mathbb{Z}.$
On the other hand, for every $q\in \mathbb{P},$ the $\mathbb{Z}$-submodule $%
K:=qG$ is coprime in $G.$ Consequently,%
\begin{equation*}
\mathrm{Spec}^{\mathrm{c}}(G)=\{qG\mid q\in \mathbb{P}\}\text{ and }\mathrm{%
Rad}^{\mathrm{c}}(G)=\bigcap\limits_{q\in \mathbb{P}}qG=0.
\end{equation*}%
Since $G$ is semisimple, every subgroup $L\leq _{\mathbb{Z}}G$ is of the
form $L=\dsum\limits_{p\in \Lambda }\mathbb{Z}/p\mathbb{Z}$ for some $%
\Lambda \subseteq \mathbb{P}.$ Let $qG\in \mathrm{Spec}^{\mathrm{c}}(G)$ be
arbitrary and suppose there are subsets $\Lambda ,\widetilde{\Lambda }%
\subseteq \mathbb{P}$ such that $\left( \dsum_{p\in \Lambda }\mathbb{Z}/p%
\mathbb{Z}\right) \cap \left( \dsum_{p\in \widetilde{\Lambda }}\mathbb{Z}/p%
\mathbb{Z}\right) \subseteq qG.$ Then indeed $q\in \Lambda \cup \widetilde{%
\Lambda }$ and it follows that $\dsum_{p\in \Lambda }\mathbb{Z}/p\mathbb{Z}%
\subseteq qG$ or $\dsum_{p\in \widetilde{\Lambda }}\mathbb{Z}/p\mathbb{Z}%
\subseteq qG.$ We conclude that every coprime subgroup of $G$ is strongly
irreducible, i.e. $\mathrm{Spec}^{\mathrm{c}}(_{\mathbb{Z}}G)\subseteq
\mathcal{SI}(_{\mathbb{Z}}G),$ whence $_{\mathbb{Z}}G$ is a top$^{\mathrm{c}}
$-module by Theorem \ref{c-top}. Notice that $_{\mathbb{Z}}G$ is not a
multiplication module: Let $\Lambda \subseteq \mathbb{P}$ be a finite subset
and $N:=\bigoplus\limits_{p\in \Lambda }\mathbb{Z}/p\mathbb{Z}.$ Then $(N:_{%
\mathbb{Z}}G)=\bigcap\limits_{p\in \mathbb{P}\backslash \Lambda }p\mathbb{Z}%
=0,$ whence $N\neq (N:_{\mathbb{Z}}G)G$ and $_{\mathbb{Z}}G$ is not a
multiplication module.
\end{ex}

\begin{lemma}
\label{c-closure}Let $_{R}M$ be a top$^{\mathrm{c}}$-module. The closure of
any $\mathcal{A}\subseteq \mathrm{Spec}^{\mathrm{c}}(M)$ is%
\begin{equation}
\overline{\mathcal{A}}=\mathcal{V}^{\mathrm{c}}(\mathcal{J}(\mathcal{A})%
\mathcal{)}.  \label{c-closure=}
\end{equation}
\end{lemma}

\begin{remarks}
\label{max-char}Let $_{R}M$ be a top$^{\mathrm{c}}$-module and consider the
Zariski topology
\begin{equation}
\mathbf{Z}^{\mathrm{c}}(M):=(\mathrm{Spec}^{\mathrm{c}}(M),\tau ^{\mathrm{c}%
}(M)).
\end{equation}

\begin{enumerate}
\item $\mathbf{Z}^{\mathrm{c}}(M)$ is a $T_{0}$ (Kolmogorov) space.

\item Set $\mathcal{X}_{m}^{\mathrm{c}}:=\mathcal{X}^{\mathrm{c}}(Rm)$ for
each $m\in M.$ The set
\begin{equation*}
\mathcal{B}:=\{\mathcal{X}_{m}^{\mathrm{c}}\mid m\in M\}
\end{equation*}%
is a basis of open sets for the Zariski topology $\mathbf{Z}^{\mathrm{c}%
}(M). $

\item If $L\in \mathrm{Spec}^{\mathrm{c}}(M),$ then $\overline{\{L\}}=%
\mathcal{V}^{\mathrm{c}}(L).$ In particular, for any $K\in \mathrm{Spec}^{%
\mathrm{c}}(M):$%
\begin{equation*}
K\in \overline{\{L\}}\Leftrightarrow L\subseteq K.
\end{equation*}

\item ${\mathcal{X}}^{\mathrm{c}}(L)=\emptyset \Rightarrow L\subseteq
\mathrm{Rad}(M).$ The converse holds if, for example, $\mathrm{Max}(M)=%
\mathrm{Spec}^{\mathrm{c}}(M).$

\item Let $_{R}M$ be coatomic. For every $L\leq _{R}M$ we have ${\mathcal{V}}%
^{\mathrm{c}}(L)=\emptyset $ if and only if $L=M.$

\item Let $L\lvertneqq _{R}M.$ The mapping%
\begin{equation*}
\pi :\mathrm{Spec}^{\mathrm{c}}(M)\rightarrow \mathrm{Spec}^{\mathrm{c}%
}(M/L),\text{ }N\mapsto (N+L)/L
\end{equation*}%
induces a continuous mapping%
\begin{equation}
\mathbf{\pi }:\mathbf{Z}^{\mathrm{c}}(M)\longrightarrow \mathbf{Z}^{\mathrm{c%
}}(M/L).
\end{equation}%
This follows from the fact that $\mathbf{\pi }^{-1}(\mathcal{V}^{\mathrm{c}%
}(N/L))=\mathcal{V}^{\mathrm{c}}(N)$ for every $L\leq _{R}N\leq _{R}M.$
\end{enumerate}
\end{remarks}

\begin{notation}
Set
\begin{equation*}
\mathbf{CL}(\mathbf{Z}^{\mathrm{c}}(M)):=\{\mathcal{A}\subseteq \mathrm{Spec}%
^{\mathrm{c}}(M)\mid \mathcal{A}=\overline{\mathcal{A}}\}\text{ and }%
\mathcal{R}^{\mathrm{c}}(M):=\{L\leq _{R}M\mid \mathrm{Rad}_{M}^{\mathrm{c}%
}(L)=L\}.
\end{equation*}
\end{notation}

\begin{theorem}
\label{c-11}Let $_{R}M$ be a top$^{\mathrm{c}}$-module.

\begin{enumerate}
\item We have an order-reversing bijection%
\begin{equation}
\mathcal{R}^{\mathrm{c}}(M)\longleftrightarrow \mathbf{CL}(\mathbf{Z}^{%
\mathrm{c}}(M)),\text{ }L\mapsto \mathcal{V}^{\mathrm{c}}(L).  \label{c-bij}
\end{equation}

\item $\mathbf{Z}^{\mathrm{c}}(M)$ is Noetherian if and only if $_{R}M$
satisfies the ACC condition on $\mathcal{R}^{\mathrm{c}}(M).$

\item $\mathbf{Z}^{\mathrm{c}}(M)$ is Artinian if and only if $_{R}M$
satisfies the DCC condition on $\mathcal{R}^{\mathrm{c}}(M).$
\end{enumerate}
\end{theorem}

\begin{theorem}
\label{c-noth-art} Let $_{R}M$ be a top$^{\mathrm{c}}$-module. If $_{R}M$ is
Noetherian \emph{(Artinian)}, then $\mathbf{Z}^{\mathrm{c}}(M)$ is
Noetherian \emph{(}Artinian\emph{)}.
\end{theorem}

\begin{punto}
Recall that $_{R}M$ is said to be \textit{distributive} iff
\begin{equation}
L\cap (K_{1}+K_{2})=(L\cap K_{1})+(L\cap K_{2})\text{ (equivalently,}%
(L+K_{1})\cap (L+K_{2})=L+(K_{1}\cap K_{2}))
\end{equation}%
for all submodules of $M.$ We call $_{R}M$ \textit{completely distributive}
iff for all $L,K_{\lambda }\in \mathcal{L}(M):$%
\begin{equation*}
\bigcap\limits_{\lambda \in \Lambda }(L+K_{\lambda
})=L+(\bigcap\limits_{\lambda \in \Lambda }K_{\lambda }).
\end{equation*}
\end{punto}

\begin{proposition}
\label{c-irr}Let $_{R}M$ be a completely distributive top$^{c}$-module and ${%
\mathcal{A}}\subseteq \mathrm{Spec}^{\mathrm{c}}(M).$

\begin{enumerate}
\item If ${\mathcal{A}}$ is irreducible, then {$\mathcal{J}$}$({\mathcal{A}}%
) $ is a coprime submodule of $M.$

\item If $\mathrm{Spec}^{\mathrm{c}}(M)\subseteq \mathcal{SI}(M),$ then the
following are equivalent:

\begin{enumerate}
\item ${\mathcal{A}}\subseteq \mathrm{Spec}^{\mathrm{c}}(M)$ is irreducible;

\item {$\mathcal{J}$}$({\mathcal{A}})$ is a coprime submodule of $M;$

\item {$\mathcal{J}$}$({\mathcal{A}})\lvertneqq _{R}M$ is strongly
irreducible.
\end{enumerate}
\end{enumerate}
\end{proposition}

\begin{Beweis}
Let $_{R}M$ be a top$^{c}$-module and ${\mathcal{A}}\subseteq \mathrm{Spec}^{%
\mathrm{c}}(M).$

\begin{enumerate}
\item Assume that ${\mathcal{A}}$ is irreducible. By definition, ${\mathcal{A%
}}\neq \emptyset $ and so {$\mathcal{J}$}$({\mathcal{A}})\lvertneqq _{R}M.$
Let $I\in \mathcal{I}(R)$ and suppose that $IM\nsubseteqq \mathcal{J}({%
\mathcal{A}})$ and $IM+\mathcal{J}({\mathcal{A}})\neq M.$ Let $\mathcal{A}%
_{1}:=\{K\in \mathcal{A}\mid IM\subseteq K\}$ and $\mathcal{A}_{2}:=\{K\in
\mathcal{A}\mid IM+K=M\}.$ Notice that ${\mathcal{A}}\subseteq \mathcal{V}^{%
\mathrm{c}}(${$\mathcal{J}(\mathcal{A}_{1}$}$))\cup \mathcal{V}^{\mathrm{c}%
}( ${$\mathcal{J}$}${(}\mathcal{A}{_{2}})).$ However, ${\mathcal{A}}%
\nsubseteqq \mathcal{V}^{\mathrm{c}}(${$\mathcal{J}(\mathcal{A}_{1}$}$))$
(otherwise, $IM\subseteq \mathcal{J}${$(\mathcal{A}_{1}$}$)=\mathcal{J}${$(%
\mathcal{A}$}$) $)\ and ${\mathcal{A}}\nsubseteqq \mathcal{V}^{\mathrm{c}}(${%
$\mathcal{J}(\mathcal{A}_{2}$}$))$ (otherwise, $IM+\mathcal{J}({\mathcal{A}}%
)=IM+\mathcal{J}({\mathcal{A}}_{2})=IM+\bigcap\limits_{K\in \mathcal{A}%
_{2}}K=\bigcap\limits_{K\in \mathcal{A}_{2}}(IM+K)=M$). This is a
contradiction, whence $IM\subseteq \mathcal{J}({\mathcal{A}})$ or $IM+%
\mathcal{J}({\mathcal{A}})=M.$ Consequently, $\mathcal{J}({\mathcal{A}})$ is
a coprime submodule of $M.$

\item $(a)\Rightarrow (b)$ follows by \textquotedblleft 1\textquotedblright .

$(b)\Rightarrow (c)$ is obvious.

$(c)\Rightarrow (a)$ Since $\mathcal{J}(\mathcal{A})$ is a proper submodule
of $M,$ we conclude that $\mathcal{A}\neq \varnothing .$ Suppose that $%
\mathcal{A}\subseteq \mathcal{V}^{\mathrm{c}}(L_{1})\cup \mathcal{V}^{%
\mathrm{c}}(L_{2})\subseteq \mathcal{V}^{\mathrm{c}}(L_{1}\cap L_{2})$ for
some $L_{1},L_{2}\leq _{R}M,$ so that $L_{2}\cap L_{2}\subseteq {\mathcal{J}}%
(\mathcal{A}).$ Since ${\mathcal{J}}(\mathcal{A})\leq _{R}M$ is strongly
irreducible, $L_{1}\subseteq {\mathcal{J}}(\mathcal{A})$ so that $\mathcal{A}%
\subseteq \mathcal{V}^{\mathrm{c}}(L_{1})$ or $L_{2}\subseteq {\mathcal{J}}(%
\mathcal{A})$ so that $\mathcal{A}\subseteq \mathcal{V}^{\mathrm{c}}(L_{2})$%
. We conclude that $\mathcal{A}$ is irreducible.$\blacksquare $
\end{enumerate}
\end{Beweis}

\begin{theorem}
\label{corad-c}Let $_{R}M$ be a completely distributive top$^{\mathrm{c}}$%
-module.

\begin{enumerate}
\item If $\mathrm{Spec}^{\mathrm{c}}(M)$ is irreducible, then $\mathrm{Rad}%
_{M}^{\mathrm{c}}(M)$ is a coprime submodule of $M.$

\item If $\mathrm{Max}(M)\subseteq \mathrm{Spec}^{\mathrm{c}}(M)$ is
irreducible, then $\mathrm{Rad}_{M}^{\mathrm{c}}(M)$ is a coprime submodule
of $M.$

\item Let $\mathrm{Spec}^{\mathrm{c}}(M)\subseteq \mathcal{SI}(M).$

\begin{enumerate}
\item The following are equivalent:

\begin{enumerate}
\item $\mathrm{Spec}^{\mathrm{c}}(M)$ is irreducible;

\item $\mathrm{Rad}_{M}^{\mathrm{c}}(M)$ is a coprime submodule of $M;$

\item $\mathrm{Rad}_{M}^{\mathrm{c}}(M)\lvertneqq _{R}M$ is strongly
irreducible.
\end{enumerate}

\item The following are equivalent:

\begin{enumerate}
\item $\mathrm{Max}(M)\subseteq \mathrm{Spec}^{\mathrm{c}}(M)$ is
irreducible;

\item $\mathrm{Rad}(M)$ is a coprime submodule of $M;$

\item $\mathrm{Rad}(M)\lvertneqq _{R}M$ is strongly irreducible.
\end{enumerate}
\end{enumerate}
\end{enumerate}
\end{theorem}

\begin{ex}
Let $_{R}M$ be a top$^{\mathrm{c}}$-module. If $\varnothing \neq {\mathcal{A}%
}\subseteq \mathrm{Spec}^{\mathrm{c}}(M)$ is a chain, then ${\mathcal{A}}$
is irreducible. In particular, if $_{R}M$ is uniserial, then\textrm{\ }$%
\mathrm{Spec}^{\mathrm{c}}(M)$ is irreducible.
\end{ex}

\begin{notation}
We denote by $\mathrm{Min}(\mathrm{Spec}^{\mathrm{c}}(M))$ the minimal
elements of $\mathrm{Spec}^{\mathrm{c}}(M).$
\end{notation}

\begin{proposition}
\label{c-max-irr}Let $_{R}M$ be completely distributive and $\mathrm{Spec}^{%
\mathrm{c}}(M)\subseteq \mathcal{SI}(M)$ (whence $_{R}M$ is a top$^{\mathrm{c%
}}$-module).

\begin{enumerate}
\item The bijection \emph{(}\ref{c-bij}\emph{)} restricts to a bijection%
\begin{equation*}
\mathrm{Spec}^{\mathrm{c}}(M)\longleftrightarrow \{{\mathcal{A}}\mid {%
\mathcal{A}}\subseteq \mathrm{Spec}^{\mathrm{c}}(M)\text{ is an irreducible
closed subset}\}.
\end{equation*}

\item The bijection \emph{(}\ref{c-bij}\emph{)} restricts to a bijection%
\begin{equation*}
\mathrm{Min}(\mathrm{Spec}^{\mathrm{c}}(M))\longleftrightarrow \{{\mathcal{A}%
}\mid {\mathcal{A}}\subseteq \mathrm{Spec}^{\mathrm{c}}(M)\text{ is an
irreducible component}\}.
\end{equation*}
\end{enumerate}
\end{proposition}

\begin{Beweis}
Recall the bijection $\mathcal{R}^{\mathrm{c}}(M)\overset{\mathcal{V}^{%
\mathrm{c}}(-)}{\longrightarrow }\mathbf{CL}(\mathbf{Z}^{\mathrm{c}}(M)).$

\begin{enumerate}
\item Let $K\in \mathrm{Spec}^{\mathrm{c}}(M).$ Then $K=\mathcal{J}(\mathcal{%
V}^{\mathrm{c}}(K))$ and it follows that the closed set $\mathcal{V}^{%
\mathrm{c}}(K)$ is irreducible by Proposition \ref{corad-c}. On the other
hand, let ${\mathcal{A}}\subseteq \mathrm{Spec}^{\mathrm{c}}(M)$ be a closed
irreducible subset. Then ${\mathcal{A}}=\mathcal{V}^{\mathrm{c}}(L)$ for
some $L\leq _{R}M.$ Notice that {$\mathcal{J}$}$({\mathcal{A}})$ is coprime
in $M$ by Proposition \ref{corad-c} and that ${\mathcal{A}}=\overline{{%
\mathcal{A}}}=\mathcal{V}^{\mathrm{c}}(${$\mathcal{J}$}$({\mathcal{A}})).$

\item Let $K$ be minimal in $\mathrm{Spec}^{\mathrm{c}}(M).$ Then $\mathcal{V%
}^{\mathrm{c}}(K)$ is an irreducible subset of $\mathrm{Spec}^{\mathrm{c}%
}(M) $ by \textquotedblleft 1\textquotedblright . By \cite{Bou1998}, $%
\mathcal{V}^{\mathrm{c}}(K)$ is contained in some irreducible component $%
\mathcal{Y}$ of $\mathrm{Spec}^{\mathrm{c}}(M).$ Since $\mathcal{Y}$ is
closed, there exists by \textquotedblleft 1\textquotedblright\ some $L\in
\mathrm{Spec}^{\mathrm{c}}(M)$ such that $\mathcal{Y}=\mathcal{V}^{\mathrm{c}%
}(L).$ If $\mathcal{V}^{\mathrm{c}}(K)\subsetneqq \mathcal{V}^{\mathrm{c}%
}(L),$ then $L\subsetneqq K,$ a contradiction. Consequently, $\mathcal{V}^{%
\mathrm{c}}(K)=\mathcal{V}^{\mathrm{c}}(L)$ is an irreducible component of $%
L.$

On the other hand, let $\mathcal{Y}$ be an irreducible component of $\mathrm{%
Spec}^{\mathrm{c}}(M).$ Then $\mathcal{Y}$ is closed and irreducible, i.e. $%
\mathcal{Y}=\mathcal{V}^{\mathrm{c}}(L)$ for some $L\in \mathrm{Spec}^{%
\mathrm{c}}(M)$ by \textquotedblleft 1\textquotedblright . Suppose that $L$
is not minimal in $\mathrm{Spec}^{\mathrm{c}}(M),$ so that there exists $%
K\in \mathrm{Spec}^{\mathrm{c}}(M)$ such that $K\subsetneqq L.$ It follows
that $\mathcal{V}^{\mathrm{c}}(L)\subsetneqq \mathcal{V}^{\mathrm{c}}(K),$ a
contradiction since $\mathcal{V}^{\mathrm{c}}(K)\subseteq \mathrm{Spec}^{%
\mathrm{c}}(M)$ is irreducible by \textquotedblleft 1\textquotedblright . We
conclude that $L$ is minimal in $\mathrm{Spec}^{\mathrm{c}}(M).\blacksquare $
\end{enumerate}
\end{Beweis}

\begin{corollary}
If $_{R}M$ is completely distributive and $\mathrm{Spec}^{\mathrm{c}%
}(M)\subseteq \mathcal{SI}(M),$ then $\mathrm{Spec}^{\mathrm{c}}(M)$ is a
Sober space.
\end{corollary}

\begin{theorem}
\label{c-hollow}Let $_{R}M$ be a coatomic top$^{\mathrm{c}}$-module. Then $%
_{R}M$ is hollow if and only if $\mathrm{Spec}^{\mathrm{c}}(M)$ is
ultraconnected.
\end{theorem}

\begin{ex}
Let $_{R}M$ be a coatomic multiplication module. Then $\mathrm{Spec}^{%
\mathrm{c}}(M)$ is ultraconnected if and only if $_{R}M$ is local. Indeed,
if $\mathrm{Spec}^{\mathrm{c}}(M)=\mathrm{Max}(M)$ is ultraconnected, then $%
\left\vert \mathrm{Max}(M)\right\vert =1:$\ If $\mathfrak{m}_{1},\mathfrak{m}%
_{2}\in \mathrm{Max}(M)$ are distinct, then $\mathcal{V}(\mathfrak{m}%
_{1})\cap \mathcal{V}(\mathfrak{m}_{2})=\{\mathfrak{m}_{1}\}\cap \{\mathfrak{%
m}_{2}\}=\varnothing ,$ a contradiction. Since $_{R}M$ is coatomic, we
conclude that $_{R}M$ is local. On the other hand, if $_{R}M$ is local, then
indeed $_{R}M$ is hollow whence $\mathrm{Spec}^{\mathrm{c}}(M)$ is
ultraconnected by Theorem \ref{c-hollow}. $\blacksquare $
\end{ex}

\begin{theorem}
\label{c-compact} Let $_{R}M$ be a coatomic top$^{\mathrm{c}}$-module.

\begin{enumerate}
\item If $\mathrm{Max}(M)$ is countable, then $\mathrm{Spec}^{\mathrm{c}}(M)$
is countably compact.

\item If $\mathrm{Max}(M)$ is finite, then $\mathrm{Spec}^{\mathrm{c}}(M)$
is compact.
\end{enumerate}
\end{theorem}

\begin{ex}
Let $M=\mathbb{Z}_{p^{\infty }}$ and recall that%
\begin{equation*}
\mathrm{Spec}^{\mathrm{c}}(\mathbb{Z}_{p^{\infty }})=\{\mathbb{Z}(\frac{1}{%
p^{n}}+\mathbb{Z})\mid n\in \mathbb{N}\}.
\end{equation*}%
Notice that $\mathrm{Max}(M)=\varnothing $ (in particular finite) but $%
\mathrm{Spec}^{\mathrm{c}}(\mathbb{Z}_{p^{\infty }})$ is not compact since
the open cover%
\begin{equation*}
\{\mathcal{X}^{\mathrm{c}}\left( \mathbb{Z}(\frac{1}{p^{n+1}}+\mathbb{Z}%
)\right) \mid n\in \mathbb{N}\}
\end{equation*}%
has no finite subcover. Recall that $\mathbb{Z}_{p^{\infty }}$ is not
coatomic as a $Z$-module. This shows that the assumption that $_{R}M$ be a
coatomic in Theorem \ref{c-compact} \textquotedblleft 2\textquotedblright\
cannot be removed.
\end{ex}

\begin{proposition}
\label{irr-c}Let $_{R}M$ be a top$^{\mathrm{c}}$-module and assume that
every coprime submodule of $M$ is maximal.

\begin{enumerate}
\item If $_{R}M$ has the complete max-property, then $\mathrm{Spec}^{\mathrm{%
c}}(M)$ is discrete.

\item $M$ has a unique maximal $R$-submodule if and only if $_{R}M$ has the
complete max-property and $\mathrm{Spec}^{\mathrm{c}}(M)$ is connected.
\end{enumerate}
\end{proposition}

\begin{theorem}
\label{c-countable}Let $_{R}M$ be a coatomic top$^{\mathrm{c}}$-module and
assume that every coprime submodule of $M$ is maximal. If $_{R}M$ has the
complete max property, then

\begin{enumerate}
\item $\mathrm{Spec}^{\mathrm{c}}(M)$ is countably compact if and only if $%
\mathrm{Max}(M)$ is countable.

\item $\mathrm{Spec}^{\mathrm{c}}(M)$ is compact if and only if $\mathrm{Max}%
(M)$ is finite.
\end{enumerate}
\end{theorem}

\qquad As a direct consequence of Proposition \ref{irr-c} we obtain:

\begin{theorem}
\label{c-colocal}Let $_{R}M$ be a coatomic top$^{\mathrm{c}}$-module with
the complete max-property and assume that every coprime submodule of $M$ is
maximal. Then $_{R}M$ is local if and only if $\mathrm{Spec}^{\mathrm{c}}(M)$
is connected.
\end{theorem}

\begin{lemma}
\label{c-1n}Let $\mathrm{Spec}^{\mathrm{c}}(M)\subseteq \mathcal{SI}(M).$ If
$n\geq 2$ and $\mathcal{A}=\{K_{1},...,K_{n}\}\subseteq \mathrm{Spec}^{%
\mathrm{c}}(M)$ is a connected subset, then for every $i\in \{1,...,n\},$
there exists $j\in \{1,...,n\}\backslash \{i\}$ such that $K_{i}\leq
_{R}K_{j}$ or $K_{j}\leq _{R}K_{i}.$
\end{lemma}

\begin{proposition}
\label{c-lf}Let $_{R}M$ be a coatomic top$^{\mathrm{c}}$-module with the
complete max-property and let $\varnothing \neq \mathcal{A}=\{K_{\lambda
}\}_{\Lambda }\subseteq \mathrm{Max}(M).$ If $\left\vert \mathcal{M}%
(L)\right\vert <\infty $ for every $L\in \mathrm{Spec}^{\mathrm{c}}(M),$
then $\mathcal{A}$ is locally finite.
\end{proposition}

\begin{proposition}
If $_{R}M$ be a coatomic top$^{\mathrm{c}}$-module, then the following are
equivalent for any $L\leq _{R}M:$

\begin{enumerate}
\item $L\in \mathrm{Max}(M);$

\item $L$ is a coprime submodule of $M$ and ${\mathcal{V}}^{\mathrm{c}%
}(L)=\{L\};$

\item $\{L\}$ is closed in ${\mathbf{Z}}_{M}^{\mathrm{c}}.$
\end{enumerate}
\end{proposition}

\begin{proposition}
\label{c-T1}Let $_{R}M$ be a coatomic top$^{\mathrm{c}}$-module. Then $%
\mathrm{Spec}^{\mathrm{c}}(M)=\mathrm{Max}(M)$ if and only if $\mathbf{Z}^{%
\mathrm{c}}(M)$ is $T_{1}$ \emph{(}Fr\'{e}cht space\emph{)}.
\end{proposition}

\qquad Combining the previous results we obtain

\begin{theorem}
\label{c-T2}Let $_{R}M$ be a coatomic top$^{\mathrm{c}}$-module with the
complete max-property. The following are equivalent:

\begin{enumerate}
\item $\mathrm{Spec}^{\mathrm{c}}(M)=\mathcal{S}(M);$

\item $\mathbf{Z}^{\mathrm{c}}(M)$ is discrete;

\item $\mathbf{Z}^{\mathrm{c}}(M)$ is $T_{2}$ \emph{(}Hausdorff space\emph{)}%
;

\item $\mathbf{Z}^{\mathrm{c}}(M)$ is $T_{1}$ \emph{(}Fr\'{e}cht space\emph{)%
}.
\end{enumerate}
\end{theorem}

\textbf{Acknowledgement}: Several results in Section 3 resulted from
discussions with Professor Patrick F. Smith during the visit of the author
to the University of Glasgow (September 2008). The author thanks Prof. Smith
for these fruitful discussions and for providing him with several examples
on coprime (coprimeless) modules. He also thanks the University of Glasgow
for the hospitality and the Deanship of Scientific Research (DSR) at King
Fahd University of Petroleum $\&$ Minerals (KFUPM) for the financial support.

\end{document}